\documentclass[10pt,twoside]{amsart}
\usepackage{amscd,amsmath,amssymb}
\usepackage[all]{xypic} 

\newtheorem{prop}{Proposition}[section] 
\newtheorem{lem}[prop]{Lemma}
\newtheorem{thm}[prop]{Theorem} 
\newtheorem{cor}[prop]{Corollary}

\theoremstyle{remark} 
\newtheorem{rem}[prop]{Remark}
\newtheorem{rems}[prop]{Remarks}
\newtheorem{ex}[prop]{Example}
\newtheorem{notation}[prop]{Notation}

\theoremstyle{definition} 
\newtheorem{defi}[prop]{Definition}

\newcommand{\ch}{{\mathcal H}}

\newcommand{\mh}{{\mathbb H}} 
\newcommand{\mq}{{\mathbb Q}}
\newcommand{\mz}{{\mathbb Z}}

\DeclareMathOperator{\disc}{disc}

\newcommand{\qform}[1]{{\langle{#1}\rangle}} 
\newcommand{\pform}[1]{{\langle\!\langle{#1}\rangle\!\rangle}} 

\newcommand{\ra}{\rightarrow}

\DeclareMathOperator{\br}{Br}
\DeclareMathOperator{\Br}{Br} 
\DeclareMathOperator{\Dec}{Dec}
\DeclareMathOperator{\cores}{cor} 
\DeclareMathOperator{\End}{End}
\DeclareMathOperator{\gal}{Gal}

\DeclareMathOperator{\ind}{ind} 
\DeclareMathOperator{\Int}{Int}
\DeclareMathOperator{\ort}{{O}}

\DeclareMathOperator{\Sym}{Sym} 
\DeclareMathOperator{\SSym}{SSym}
\DeclareMathOperator{\so}{O^+}
\DeclareMathOperator{\Nrd}{Nrd}
\newcommand{\sq}[1]{#1^\times/#1^{\times 2}}

\DeclareMathOperator{\GL}{GL}

\newcommand{\Spin}{\mathrm {Spin}} 

\DeclareMathOperator{\SL}{SL}

\DeclareMathOperator{\ad}{ad} 
\DeclareMathOperator{\Ad}{Ad}
\newcommand{\ba}{\overline{\rule{2.5mm}{0mm}\rule{0mm}{4pt}}} 

\newcommand{\ot}{\otimes}
\newcommand{\QZt}{\mathbb{Q}/\mathbb{Z}(2)} 

\newcommand{\hyp}{{\mathrm {hyp}}}

\newcommand{\Fields}{\operatorname{\textsf{Fields}}}

\def\orthsum{\operatornamewithlimits{%
  \mathchoice{\boxplus}{\boxplus}{\boxplus}{\boxplus}}}

\def\YEAR{\year}\newcount\VOL\VOL=\YEAR\advance\VOL by-1995
\def\firstpage{1}\def\lastpage{46}
\def\received{June 30, 2014}\def\revised{}
\def\communicated{}

\makeatletter
\def\magnification{\afterassignment\m@g\count@}
\def\m@g{\mag=\count@\hsize6.5truein\vsize8.9truein\dimen\footins8truein}
\makeatother

\font\eightrm=cmr8
\font\caps=cmcsc10                    
\font\Caps=cmcsc10 scaled \magstep1   

\makeatletter
\@specialpagefalse\headheight=8.5pt
\def\DocMath{}

\renewcommand{\@evenhead}{%
    \ifnum\thepage>\lastpage\rlap{\thepage}\hfill%
    \else\rlap{\thepage}\slshape\leftmark\hfill{\caps\SAuthor}\hfill\fi}%
\renewcommand{\@oddhead}{%
    \ifnum\thepage=\firstpage{\DocMath\hfill\llap{\thepage}}%
    \else{\slshape\rightmark}\hfill{\caps\STitle}\hfill\llap{\thepage}\fi}%
\makeatother

\def\TSkip{\bigskip}
\newbox\TheTitle{\obeylines\gdef\GetTitle #1
\ShortTitle  #2
\SubTitle    #3
\Author      #4
\ShortAuthor #5
\EndTitle
{\setbox\TheTitle=\vbox{\baselineskip=20pt\let\par=\cr\obeylines%
\halign{\centerline{\Caps##}\cr\noalign{\medskip}\cr#1\cr}}%
	\copy\TheTitle\TSkip\TSkip%
\def\next{#2}\ifx\next\empty\gdef\STitle{#1}\else\gdef\STitle{#2}\fi%
\def\next{#3}\ifx\next\empty%
    \else\setbox\TheTitle=\vbox{\baselineskip=20pt\let\par=\cr\obeylines%
    \halign{\centerline{\caps##} #3\cr}}\copy\TheTitle\TSkip\TSkip\fi%
\setbox\TheTitle=\vbox{\let\par=\cr\obeylines%
\halign{\centerline{\caps##} #4\cr}}\copy\TheTitle\TSkip\TSkip%
\def\next{#5}\ifx\next\empty\gdef\SAuthor{#4}\else\gdef\SAuthor{#5}\fi%
\ifx\received\empty\relax
    \else\centerline{\eightrm Received: \received}\fi%
\ifx\revised\empty\TSkip%
    \else\centerline{\eightrm Revised: \revised}\TSkip\fi%
\ifx\communicated\empty\relax
    \else\centerline{\eightrm Communicated by \communicated}\fi\TSkip\TSkip%
\catcode'015=5}}\def\Title{\obeylines\GetTitle}
\def\Abstract{\begingroup\narrower
    \parskip=\medskipamount\parindent=0pt{\caps Abstract. }}
\def\EndAbstract{\par\endgroup\TSkip}

\long\def\MSC#1\EndMSC{\def\arg{#1}\ifx\arg\empty\relax\else
     {\par\narrower\noindent%
     2010 Mathematics Subject Classification: #1\par}\fi}

\long\def\KEY#1\EndKEY{\def\arg{#1}\ifx\arg\empty\relax\else
	{\par\narrower\noindent Keywords and Phrases: #1\par}\fi\TSkip}

\newbox\TheAdd\def\Addresses{\vfill\copy\TheAdd\vfill
    \ifodd\number\lastpage\vfill\eject\phantom{.}\vfill\eject\fi}
{\obeylines\gdef\GetAddress #1
\Address #2 
\Address #3
\Address #4
\EndAddress
{\def\xs{5.3truecm}\parindent=0pt
\setbox0=\vtop{{\obeylines\hsize=\xs#1\par}}\def\next{#2}
\ifx\next\empty 
     \setbox\TheAdd=\hbox to\hsize{\hfill\copy0\hfill}
\else\setbox1=\vtop{{\obeylines\hsize=\xs#2\par}}\def\next{#3}
\ifx\next\empty 
     \setbox\TheAdd=\hbox to\hsize{\hfill\copy0\hfill\copy1\hfill}
\else\setbox2=\vtop{{\obeylines\hsize=\xs#3\par}}\def\next{#4}
\ifx\next\empty\ 
     \setbox\TheAdd=\vtop{\hbox to\hsize{\hfill\copy0\hfill\copy1\hfill}
                \vskip20pt\hbox to\hsize{\hfill\copy2\hfill}}
\else\setbox3=\vtop{{\obeylines\hsize=\xs#4\par}}
     \setbox\TheAdd=\vtop{\hbox to\hsize{\hfill\copy0\hfill\copy1\hfill}
	        \vskip20pt\hbox to\hsize{\hfill\copy2\hfill\copy3\hfill}}
\fi\fi\fi\catcode'015=5}}\gdef\Address{\obeylines\GetAddress}

\begin{document}
\Title
The Arason invariant of orthogonal involutions of degree 12 and 8, 
and quaternionic subgroups of the Brauer group
\ShortTitle 
Arason invariant and quaternionic subgroups
\SubTitle   
To Sasha Merkurjev on his 60th birthday
\Author 
Anne Qu\'eguiner-Mathieu$^1$ and Jean-Pierre Tignol$^2$
\ShortAuthor 
A. Qu\'eguiner-Mathieu, J.-P. Tignol
\EndTitle
\Abstract 
  Using the Rost invariant for torsors under $\Spin$ groups one may
  define an analogue of the Arason invariant for certain hermitian
  forms and orthogonal involutions. We calculate this invariant explicitly in various cases, and
  use it to associate to every orthogonal involution $\sigma$ with
  trivial discriminant and trivial Clifford invariant over a central
  simple algebra $A$ of even co-index an element $f_3(\sigma)$ in the
  subgroup $F^\times\cdot[A]$ of $H^3(F,\QZt)$. This invariant $f_3(\sigma)$ is
  the double of any representative of the Arason invariant $e_3(\sigma)\in H^3(F,\QZt)/F^\times\cdot[A]$; it vanishes when $\deg A\leq10$ and also when
  there is a quadratic extension of $F$ that simultaneously splits $A$
  and makes $\sigma$ hyperbolic.  
  The paper provides a detailed study of both invariants, with particular attention to the degree $12$ case, and to the relation with the existence of a quadratic splitting field.

   As a main tool 
  we establish, when $\deg(A)=12$, an additive decomposition of $(A,\sigma)$ into three
  summands that are central simple algebras of degree~$4$ with
  orthogonal involutions with trivial discriminant, extending a
  well-known result of Pfister on quadratic forms of dimension~$12$ in
  $I^3F$. The Clifford components of the summands generate a subgroup $U$ of the Brauer group of $F$, in which every
  element is represented by a quaternion algebra. We show that the Arason invariant 
  $e_3(\sigma)$ generates the homology of a complex of degree $3$ Galois cohomology groups, attached to the subgroup $U$, 
  which was introduced and studied by Peyre. 
  In the final section, we use the results
  on degree~$12$ algebras to extend the definition of the Arason invariant to trialitarian triples in which all three algebras have index at most $2$.

\EndAbstract
\MSC 
11E72, 11E81, 16W10.
\EndMSC
\KEY 
Cohomological invariant, orthogonal group, algebra with involution,
Clifford algebra.
\EndKEY
\Address 
Anne Qu\'eguiner-Mathieu
Universit\'e Paris 13, 
Sorbonne Paris Cit\'e
LAGA - CNRS (UMR 7539) 
F-93430 Villetaneuse, France
queguin@math.univ-paris13.fr
\Address
Jean-Pierre Tignol
ICTEAM Institute, Box L4.05.01
Universit\'e catholique de Louvain
B-1348 Louvain-la-Neuve, Belgium
jean-pierre.tignol@uclouvain.be
\Address

\Address

\EndAddress
\footnotetext[1]{The first author acknowledges the support of the French Agence Nationale de la Recherche (ANR) under reference ANR-12-BL01-0005.}
\footnotetext[2]{The second author is grateful to the first author
  and the Universit\'e Paris 13 for their hospitality while the work
  for this paper was carried out. He acknowledges support
    from the Fonds de la Recherche Scientifique--FNRS under grants
    n$^\circ$~1.5009.11 and 1.5054.12.}


\section{Introduction}

In quadratic form theory, the Arason invariant is a 
degree~$3$ Galois cohomology class with $\mu_2$ coefficients
attached to an even-dimensional quadratic form with trivial
discriminant and trivial Clifford invariant. Originally defined by
Arason in~\cite{Arason}, it can also be described in terms of the Rost
invariant of a split Spin group, as explained
in~\cite[\S\,31.B]{KMRT}.  
It is not always possible to extend this invariant to the more general
setting of orthogonal involutions, see~\cite[\S3.4]{BFPQM}.
Nevertheless, one may use the Rost invariant of some possibly
non-split Spin groups to define {\em relative} and {\em absolute}
Arason invariants for some orthogonal involutions (see~\cite{Hyd} or
section~\ref{section:Arason} below for precise definitions).  
This was first noticed by Bayer-Fluckiger and Parimala
in~\cite{BFP98},
where they use the Rost invariant to prove classification theorems for
hermitian or skew-hermitian forms, leading to a proof of the
so-called Hasse Principle conjecture~II.  

In this paper, we call Arason invariant the absolute Arason invariant.
For orthogonal involutions, it was considered by Garibaldi, who uses
the notation $e_3^{\mathrm {hyp}}$, in~\cite{Ga}, and by Berhuy in the
index $2$ case in~\cite{Berhuy}. In particular, the latter covers the
case of central simple algebras of degree $2m$ with $m$ odd, since
such an algebra has index $1$ or $2$ when it is endowed with an
orthogonal involution.

Based on the Rost invariant for the exceptional group $E_8$, Garibaldi
also defined, for orthogonal involutions on degree $16$ central simple
algebras, another invariant related to the Arason invariant of
quadratic forms, denoted by $e_3^{16}$. Bermudez and Ruozzi~\cite{BR}
extended this definition to all degrees divisible by $16$. It follows
from the proof of Corollary 10.11 in~\cite{Ga}, and Remark 4.10 in
Barry's paper~\cite{Barry}, that these invariants do \emph{not}
coincide with what we call Arason invariant in this paper.

A systematic study of the relative and absolute Arason invariants for
orthogonal involutions was recently initiated in~\cite{QT:Arason},
where the degree $8$ case is studied in detail. In this paper, we
continue with an investigation of absolute invariants in degree $12$.  

Let $(A,\sigma)$ be a central simple algebra with orthogonal
involution over a field $F$ of characteristic different
  from~$2$. The Arason invariant $e_3(\sigma)$, when
defined, belongs to the
quotient $$H^3(F,\mq/\mz(2))/F^\times\cdot[A],$$ where
$F^\times\cdot[A]$ denotes the subgroup consisting of cup products
$(\lambda)\cdot [A]$, for $\lambda\in F^\times$, $[A]$ the Brauer
class of $A$. In \S\,\ref{section:Arason} below, we give a general
formula for computing the Arason invariant of an algebra with
involution admitting a rank $2$ factor. It follows from this formula
that the Arason invariant is not always represented by a cohomology
class of order $2$. This reflects the fact that the Dynkin index of a
non-split Spin group, in large enough degree, is equal to $4$. We
define a new invariant $f_3(\sigma)\in H^3(F,\mu_2)$, attached to any
orthogonal involution for which the Arason invariant is defined, and
which vanishes if and only if the Arason invariant is
represented by a cohomology class of order $2$. This invariant is zero
if the algebra is split, or of degree $\leq 10$; starting in degree
$12$, we produce explicit examples where it is non-zero. This is an
important motivation for studying the degree $12$ case in detail.  

The main results of the paper are given in sections~\ref{sec:adddec12} to~\ref{sec:quadsplit}. 
First, we prove that a degree $12$ algebra with orthogonal involution
$(A,\sigma)$, having trivial discriminant and trivial Clifford
invariant, admits a non-unique decomposition as a sum---in the sense of
algebras with involution---of three degree $4$ algebras with orthogonal
involution of trivial discriminant.  
This can be seen as a refinement of the main result
  of~\cite{GQ:deg12}, even though our proof in index $4$
relies on the open-orbit argument of \cite{GQ:deg12}, see
Remark~\ref{GaQ.rem}. Using this additive decomposition, we
  associate to $(A,\sigma)$ in a non-canonical
way some subgroups of the Brauer group of $F$, which we call
decomposition groups of $(A,\sigma)$, see
Definition~\ref{defi:decsubgrp}.  
Such subgroups $U\subset\br(F)$ are
generated by (at most) three quaternion algebras; they were considered
by Peyre in~\cite{Peyre}, where the homology of the following
complex is studied:
$$F^\times\cdot U\ra H^3(F,\mq/\mz(2))\ra H^3(F_U,\mq/\mz(2)), $$
where $F^\times\cdot U$ denotes the subgroup generated by cup products $(\lambda)\cdot [B]$, for $\lambda\in F^\times$,  $[B]\in U$, and 
$F_U$ is the function field of the product of the Severi-Brauer
varieties associated to the elements of $U$. Peyre's results
are recalled in \S\,\ref{sec:Peyre}.  

In~\S\,\ref{sec:ArasonPeyre}, we restrict to those algebras with
involution of degree $12$ for which the Arason invariant is defined,
and we prove $e_3(\sigma)$ detects isotropy of $\sigma$, and vanishes
if and only if $\sigma$ is hyperbolic.  
We then explore the relations between the decomposition groups and the
values of the Arason invariant. Reversing the viewpoint, we also prove
that the Arason invariant $e_3(\sigma)$ provides a generator of the
homology of Peyre's complex, where $(A,\sigma)$ is any algebra with
involution admitting $U$ as a decomposition group.  

In~\S\,\ref{sec:quadsplit}, we give a necessary and sufficient
condition for the vanishing of $f_3(\sigma)$ in  
degree $12$, in terms of decomposition groups of $(A,\sigma)$. 
When there is a quadratic extension that splits $A$ and makes
  $\sigma$ hyperbolic, an easy corestriction argument shows that
  $f_3(\sigma)=0$, see Proposition~\ref{quadsplit.prop}. We give in
  \S\,\ref{ex.sec} an explicit example to show that the converse does
  not hold.
This also provides new examples of subgroups $U$ for which the
homology of Peyre's complex is non-trivial, which differ from Peyre's
example in that the homology is generated by a Brauer class of order
$2$.  

In the last section, we extend the definition of the Arason invariant
in degree $8$ to index $2$ algebras with involution of trivial
discriminant, and such that the two components of the Clifford algebra
have index $2$. In this case also, the algebra with involution has an
additive decomposition, and the Arason invariant detects isotropy. 

\subsection*{Acknowledgements} 
Both authors thank A. Sivatski for pointing out a mistake in a
preliminary version of the paper, and S. Garibaldi for useful
comments.

\subsection*{Notation}

Throughout this paper, we work over a base field $F$ of characteristic
different from $2$. 
We use the notation $H^n(F,M)=H^n({\mathrm
  {Gal}}(F_{\mathrm{sep}}/F),M)$ for any discrete torsion Galois
module $M$. For every integer $n\geq0$ we let
\[
H^n(F)=H^n\bigl(F,\mq/\mz(n-1)\bigr),
\]
(see \cite[Appendix~A, p.~151]{GMS} for a precise definition). The cohomology classes we consider
actually are in the $2$-primary part of these groups, hence we shall
not need the modified definition for the $p$-primary part when
$\operatorname{char}(F)=p\neq0$. For each integer $m\geq0$ we let
${}_m H^n(F)$ denote the $m$-torsion subgroup of $H^n(F)$. Using the norm-residue isomorphism, one may check that 
\[ {}_2H^n(F)=H^n(F,\mu_2)\qquad\text{and}\qquad
{}_4H^3(F)=H^3(F,\mu_4^{\otimes 2}),
\]
(see for instance~\cite[Remark 4.1]{Peyre}). 
In particular, ${}_2H^1(F)=\sq{F}$. For every $a\in F^\times$ we let
$(a)\in{}_2H^1(F)$ be the square class of $a$. For $a_1$, \ldots,
$a_n\in F^\times$ we let $(a_1,\ldots,a_n)\in{}_2H^n(F)$ be the
cup-product
\[
(a_1,\ldots,a_n)=(a_1)\cdot\:\cdots\:\cdot(a_n).
\]

We refer to \cite{KMRT} and to \cite{Lam} for background information
on central simple algebras with involution and on quadratic
forms. However, we depart from the notation in \cite{Lam} by letting
$\pform{a_1,\ldots,a_n}$ denote the $n$-fold Pfister form
\[
\pform{a_1,\ldots,a_n}=\qform{1,-a_1}\cdot\: \cdots\: \cdot
\qform{1,-a_n} \qquad\text{for $a_1$, \ldots, $a_n\in F^\times$.}
\]
Thus, the discriminant, the Clifford invariant and the Arason
invariant, viewed as cohomological invariants $e_1$, $e_2$ and $e_3$,
satisfy:
\[
e_1\bigl(\pform{a_1}\bigr)=(a_1),\quad
e_2\bigl(\pform{a_1,a_2}\bigr)=(a_1,a_2),\quad
e_3\bigl(\pform{a_1,a_2,a_3}\bigr)=(a_1,a_2,a_3).
\]

For every central simple $F$-algebra $A$, we let $[A]$ be the Brauer
class of $A$, which we identify to an element in $H^2(F)$. If $L$ is a
field extension of $F$, we let $A_L=A\otimes_FL$ be the $L$-algebra
obtained from $A$ by extending scalars. 

Recall that the object function from the category ${\mathrm
  {Fields}}_F$ of field extensions of $F$ to abelian groups defined by
\[
L\mapsto\coprod_{n\geq 0} H^n(L)
\]
is a cycle module over ${\mathrm {Spec}} F$
(see~\cite[Rem.1.11]{Rost}). In particular, each group $H^n(L)$ is a
module over the Milnor $K$-ring $K_*L$.  The Brauer class $[A]$ of the
algebra $A$ generates a cycle submodule; we let $M_A$ denote the
quotient cycle module. Thus, for every field $L\supseteq F$, we have
\[
M^n_A(L) =
\begin{cases}
  H^n(L)&\text{if $n=0$ or $1$;}\\
  H^n(L)/\bigl(K_{n-2}L\cdot [A_L]\bigr) &\text{if $n\geq2$.}
\end{cases}
\]
In particular, $M^2_A(L) = \br(L)/\{0,[A_L]\}$ and
$M^3_A(L)=H^3(L)/\bigl(L^\times\cdot[A_L]\bigr)$.

Let $F_A$ denote the function field of the Severi--Brauer variety of
$A$, which is a generic splitting field of $A$. Scalar extension from
$F$ to $F_A$ yields group homomorphisms
\[
M^2_A(F)\to M^2_A(F_A)=\br(F_A) \qquad\text{and}\qquad M_A^3(F)\to
M^3_A(F_A)=H^3(F_A).
\]
The first map is injective by Amitsur's theorem,
see~\cite[Th.~5.4.1]{GS}; the second one is injective if the Schur
index of $A$ divides $4$ or if $A$ is a division algebra that
decomposes into a tensor product of three quaternion algebras, but it
is not always injective (see \cite{Peyre}, \cite{Karp}
and~\cite{Karp-triquat}).

\section{Cohomological invariants of orthogonal forms and involutions}
\label{section:Arason}

Most of this section recalls well-known facts on absolute and relative
Arason invariants that will be used in the sequel of the paper. Since
we will consider additive decompositions of algebras with involution,
we need to state the results both for hermitian forms and for
involutions.  Some new results are also included. In
Proposition~\ref{compute.prop} and Corollary~\ref{compute.cor}, we
give a general formula for the Arason invariant of an algebra with
involution which has a rank~$2$ factor. In
Definitions~\ref{f3.her} and~\ref{f3.inv}, we introduce a new
invariant, called the $f_3$-invariant, which detects whether the
Arason invariant is represented by a cohomology class of order $2$.
Finally, we state and prove in Proposition~\ref{prop:f3sum} a general
formula for computing the $f_3$ invariant of a sum of hermitian forms,
which is used in the proof of the main results of the paper.

Throughout this section, $D$ is a central division algebra over an
arbitrary field $F$ of characteristic different from $2$, and $\theta$
is an $F$-linear involution on $D$ (i.e., an involution of the first
kind). To any nondegenerate
hermitian or skew-hermitian module $(V,h)$ over $(D,\theta)$ we may
associate the corresponding adjoint algebra with involution
$\Ad_h=(\End_DV,\ad_h)$. Conversely, any central simple algebra $A$
over $F$ Brauer-equivalent to $D$ and endowed with an $F$-linear
involution $\sigma$ can be represented as $(A,\sigma)\simeq \Ad_h$ for
some nondegenerate hermitian or skew-hermitian module $(V,h)$ over
$(D,\theta)$. The hermitian or skew-hermitian module $(V,h)$ is said
to be a \emph{hermitian module of orthogonal type} if the
adjoint involution $\ad_h$ on $\End_DV$ is of orthogonal type. This
occurs if and only if either $h$ is hermitian and $\theta$ is of
orthogonal type, or $h$ is skew-hermitian and $\theta$ is of
symplectic type, see \cite[(4.2)]{KMRT}. Abusing terminology, we also
say 
that $h$ is a \emph{hermitian form of orthogonal type} when $(V,h)$ is
a hermitian module of orthogonal type (even though $h$ may actually be
skew-hermitian if $\theta$ is symplectic).

\subsection{Invariants of hermitian forms of orthogonal type}

Let $(V,h)$ be a hermitian module of orthogonal type over
$(D,\theta)$; we call $r=\dim_DV$ the \emph{relative rank} of $h$
and $n=\deg\End_DV$ the \emph{absolute rank} of $h$. These
invariants are related by $n=r\deg D$. Cohomological invariants of $h$
are defined in terms of invariants of the adjoint involution 
$\ad_h$. Namely, if $n$ is even, the \emph{discriminant} of $h$, denoted
$e_1(h)\in H^1(F,\mu_2)$, is the discriminant of $\ad_h$; the
corresponding quadratic \'etale extension $K/F$ is called the
\emph{discriminant extension}.  If $n$ is even and
$e_1(h)$ is trivial, the \emph{Clifford invariant} of $h$, denoted
$e_2(h)$, is the class in $M_D^2(F)$ of any component of the Clifford
algebra of $\ad_h$.

\begin{rem}
  It follows from the relations between the components of the Clifford
  algebra (see~\cite[(9.12)]{KMRT}) that the Clifford invariant is
  well-defined.  However, since we do not assume $n$ is divisible
  by $4$, this invariant need not be represented by a cohomology class
  of order $2$ in general.
\end{rem}

Our definitions of rank and discriminant differ slightly from the
definitions used by Bayer and Parimala in~\cite[\S 2]{BFP}, who call
``rank''  what we call the relative rank of $h$. The 
discriminant $d(h)$ of $h$ in the sense of~\cite[\S 2.1]{BFP} is
related to $e_1(h)$ by
$$
e_1(h)=d(h)\disc(\theta)^r,
$$
where $\disc(\theta)$ is the discriminant of $\theta$ as defined
in~\cite[\S7]{KMRT}, and $H^1(F,\mu_2)$ is identified with the
group of square classes $\sq F$.
In particular, $e_1(h)=d(h)$ when $h$ has even
relative rank $r$. 
By~\cite[2.1.3]{BFP}, the Clifford
invariant $\mathcal{C}\ell(h)$ used by Bayer and Parimala coincides with
our $e_2(h)$ when they are both defined, i.e., when $h$ has even
relative rank and trivial discriminant.
Assume now that the hermitian form $h$ has even relative rank, i.e.,
$\dim_DV$ is even.
The vector space $V$ then carries a hyperbolic hermitian form $h_0$ of
orthogonal type, and the standard nonabelian Galois cohomology technique
yields a canonical bijection between $H^1(F,\ort(h_0))$
and the set of isomorphism classes of nondegenerate hermitian forms of
orthogonal 
type on $V$, under which the trivial torsor
corresponds to the isomorphism class of $h_0$, see
\cite[\S29.D]{KMRT}. If $e_1(h)$ and $e_2(h)$ are trivial, the torsor
corresponding to the isomorphism class of $h$ has two
different lifts to $H^1(F,\ort^+(h_0))$, and one of
these lifts can be further lifted to a torsor $\xi$ in
$H^1(F,\Spin(h_0))$. Bayer and Parimala consider the
Rost invariant $R(\xi)\in H^3(F)$ and define in~\cite[\S 3.4,
p.~664]{BFP98} an \emph{Arason invariant} of $h$ by the formula
$$
e_3(h)=R(\xi)+F^\times\cdot[D]\in M^3_D(F). 
$$
This invariant satisfies the following properties:

\begin{lem}[Bayer--Parimala~{\cite[Lemma~3.7, Corollary~3.9]{BFP98}}]
  \label{Arason-herm.lem}
  Let $h$ and $h'$ be two hermitian forms of orthogonal type over
  $(D,\theta)$ with even relative rank, trivial discriminant, and
  trivial Clifford invariant.
  \begin{enumerate}
  \item[(i)] If $h$ is hyperbolic, then $e_3(h)=0$;
  \item[(ii)] $e_3(h\perp h')=e_3(h)+e_3(h')$;
  \item[(iii)] $e_3(\lambda h)=e_3(h)$ for any $\lambda\in F^\times$.
  \end{enumerate}
\end{lem}

In particular, it follows immediately that $e_3(h)$ is a well-defined
invariant of the Witt class of $h$. Moreover, we have:

\begin{cor}
  \label{cor:ord2}
  The Arason invariant $e_3$ has order $2$.
\end{cor}

\begin{proof}
  Indeed, for any $h$ as above, we have
  $2e_3(h)=e_3(h)+e_3(h)=e_3(h)+e_3(-h)=e_3(h\perp (-h))=0$, since
  $h\perp(-h)$ is hyperbolic.
\end{proof}

Using the properties of the Arason invariant, we may define a new
invariant as follows. Assume $h$ is as above, a hermitian form of
orthogonal type with
even relative rank, trivial discriminant, and trivial Clifford
invariant.  Let $c$, $c'\in H^3(F)$ be two representatives of the Arason
invariant $e_3(h)$. Since $c-c'\in F^\times\cdot [D]$, we have
$2c=2c'\in H^3(F)$, hence $2c$ depends only on $h$ and not on the
choice of the representative $c$ of $e_3(h)$. Because of
Corollary~\ref{cor:ord2}, the image of $2c$ in $M^3_D(F)$ vanishes,
hence $2c\in F^\times\cdot[D]$. These observations lead to the
following definition: 

\begin{defi} 
\label{f3.her}
Given an arbitrary representative $c\in H^3(F)$ of the Arason
invariant $e_3(h)\in M^3_D(F)$, we let $f_3(h)=2c\in
F^\times\cdot[D]\subset {}_2H^3(F)$. 
\end{defi}

Thus, the invariant $f_3(h)$ is well-defined; it vanishes if and only if
the Arason invariant $e_3(h)$ is represented by a class of order at
most $2$, or equivalently, if every representative of $e_3(h)$ is a
cohomology class of order at most $2$. It is clear from the definition
that the $f_3$ invariant is trivial when $D$ is split.
Another case where the $f_3$ invariant vanishes is the following:

\begin{prop}
  \label{quadsplit.prop}
  If there exists a quadratic extension $K/F$ such that $D_K$ is
  split and $h_K$ is hyperbolic, then $f_3(h)=0$.
\end{prop}

\begin{proof}
Assume such a field $K$ exists, and let $c\in H^3(F)$ be any
representative of $e_3(h)\in M_D^3(F)$.  Since $h_K$ is
hyperbolic, we have $e_3(h_K)=c_K=0\in M_D^3(K)=H^3(K)$. Hence,
$\cores_{K/F}(c_K)=2c=0$, that is $f_3(h)=0$.
\end{proof}

We will see in \S\ref{sec:quadsplit} that the converse of
Proposition~\ref{quadsplit.prop} does not hold, even in absolute rank
$12$, which is the smallest absolute rank where the $f_3$ invariant
can be nonzero. 

 Since the Dynkin index of the group $\Spin(\Ad_{h_0})$ divides $4$, 
  the Arason invariant $e_3(h)$ is represented by a
  cohomology class of order dividing $4$. Moreover, there are examples
  where it is represented by a cohomology class of order equal to $4$.
  Therefore, $f_3(h)$ is
  nonzero in general. Explicit examples can be constructed by means of
  Proposition~\ref{compute.prop} below, which yields the $e_3$ and $f_3$
  invariants of hermitian forms with a rank~$2$ factor. (See
  also Corollary~\ref{f3.cor} for examples in the lowest
  possible degree, which is~$12$.)  

\subsection{Hermitian forms with a rank~$2$ factor}

Consider a hermitian form which admits a decomposition as 
$\qform{1,-\lambda}\otimes h$ for some $\lambda\in F^\times$ and some hermitian form $h$. 
In this case, we have the following explicit formulae for the Arason and the $f_3$-invariant, when they are defined: 

\begin{prop}
  \label{compute.prop}
  Let $h$ be a hermitian form of orthogonal type with even absolute
  rank $n$, and let $K/F$ 
  be the discriminant quadratic extension.  For any $\mu\in K^\times$,
  the hermitian form $\qform{1,-N_{K/F}(\mu)} h$ has even relative
  rank, trivial discriminant and trivial Clifford invariant. Moreover,
  $$
  e_3\bigl(\qform{1,-N_{K/F}(\mu)}
  h\bigr)=\cores_{K/F}\bigl(\mu\cdot e_2(h_K)\bigr)
  $$
  and
  $$
  f_3\bigl(\qform{1,-N_{K/F}(\mu)}
  h\bigr)=
  \begin{cases}
    0&\text{if $n\equiv0\bmod4$,}\\
    N_{K/F}(\mu)\cdot[D]&\text{if $n\equiv2\bmod4$.}
  \end{cases}
  $$
  In particular, if $h$ has trivial discriminant, then for $\lambda\in
  F^\times$ we have
  $$
  e_3\bigr(\qform{1,-\lambda} h\bigl)=\lambda\cdot e_2(h)
  $$
  and
  $$
  f_3\bigr(\qform{1,-\lambda} h\bigl)=
  \begin{cases}
    0&\text{if $n\equiv0\bmod4$,}\\
    \lambda\cdot[D]&\text{if $n\equiv2\bmod4$.}
  \end{cases}
  $$
\end{prop}

\begin{proof}
  We first need to prove that the hermitian form
  $\qform{1,-N_{K/F}(\mu)}h$ has trivial discriminant and trivial
  Clifford invariant. This can be checked after scalar extension to a
  generic splitting field of $D$, since the corresponding restriction
  maps ${}_2H^1(F)\ra {}_2H^1(F_D)$ and $M_D^2(F)\ra H^2(F_D)$ are
  injective. In the split case, the result follows from an easy
  computation for the discriminant, and from~\cite[Ch.~V, \S3]{Lam}
  for the 
  Clifford invariant.  Alternatively, one may observe that the algebra
  with involution $\Ad_{\qform{1,-N_{K/F}(\mu)}h}$ decomposes as
  $\Ad_\qform{1,-N_{K/F}(\mu)}\otimes\Ad_h$, and
  apply~\cite[(7.3)(4)]{KMRT} and~\cite{Tao}. This computation also
  applies to the trivial discriminant case, where $\lambda=N_{F\times
    F/F}(\lambda, 1)$.
 
  With this in hand, we may compute the Arason invariant by using the
  description of cohomological invariants of quasi-trivial tori given
  in~\cite{MPT}.  Let us first assume $h$ has trivial
  discriminant. Consider the multiplicative group scheme
  $\mathbb{G}_m$ as a functor from the category ${\mathsf{Fields}}_F$
  to the category of abelian groups. For any field $L$ containing $F$,
  consider the map
  \[
  \varphi_L\colon \mathbb{G}_m(L)\to M^3_D(L) \quad\text{defined by}
  \quad
  \lambda\mapsto e_3\bigl(\qform{1,-\lambda} h_L\bigr).
  \]
  To see that $\varphi_L$ is a group homomorphism, observe that in the
  Witt group of $D_L$ we have for $\lambda_1$, $\lambda_2\in L^\times$
  \[
  \qform{1,-\lambda_1\lambda_2} h_L = \qform{1,-\lambda_1} h_L +
  \qform{\lambda_1}\qform{1,-\lambda_2} h_L.
  \]
  Therefore, Lemma~\ref{Arason-herm.lem} yields
  \[
  e_3\bigl(\qform{1,-\lambda_1\lambda_2}h_L\bigr) =
  e_3\bigl(\qform{1,-\lambda_1}h_L\bigr) +
  e_3\bigl(\qform{1,-\lambda_2}h_L\bigr).
  \]
  The collection of maps $\varphi_L$ defines a natural transformation
  of functors $\mathbb{G}_m\to M^3_D$, i.e., a degree~3 invariant of
  $\mathbb{G}_m$ with values in the cycle module $M_D$. By
  \cite[Prop.~2.5]{Merk}, there is an element $u\in M^2_D(F)$ such
  that for any $L$ and any $\lambda\in L^\times$
  \begin{equation*}
    \label{eq:ehypad1}
    \varphi_L(\lambda) = \lambda\cdot u_L
    \qquad\text{in $M_D^3(L)$.}
  \end{equation*}
  To complete the computation of $e_3\bigl(\qform{1,-\lambda}h\bigr)$,
  it only remains to show that 
  $u=e_2(h)$. Since the restriction map $M_D^2(F)\to
  M_D^2(F_D)=H^2(F_D)$ is injective, it suffices to show that
  $u_{F_D}=e_2(h)_{F_D}$. Now, since $F_D$ is a splitting field for
  $D$, there exists a quadratic form $q$ over $F_D$, with trivial
  discriminant, such that $({\Ad_h})_{F_D}\simeq\Ad_q$.  Let $t$ be an
  indeterminate over $F_D$. We have
  \[
  \Ad_\qform{1,-t}\otimes(\Ad_h)_{F_D(t)}\simeq
  \Ad_{\qform{1,-t}}\otimes (\Ad_q)_{F_D(t)}
  \]
  hence $e_3(\qform{1,-t}h_{F_D(t)})$ is the Arason invariant of the
  quadratic form $\qform{1,-t} q_{F_D(t)}$, which is $t\cdot
  e_2(q)=t\cdot e_2(h)_{F_D(t)}$. Therefore, we have
  \[
  t\cdot u_{F_D(t)} = t\cdot e_2(h)_{F_D(t)}.
  \]
  Taking the residue $\partial\colon H^3(F_D(t))\to H^2(F_D)$ for the
  $t$-adic valuation, we obtain $u_{F_D}=e_2(h)_{F_D}$, which
  completes the proof of the formula for
  $e_3\bigl(\qform{1,-\lambda}h\bigr)$. To compute
  $f_3\bigl(\qform{1,-\lambda}h\bigr)$, recall that $e_2(h)$ is
  represented by any of the two components $C_+$, $C_-$ of the
  Clifford algebra of $\Ad_h$. Therefore,
  $e_3\bigl(\qform{1,-\lambda}h\bigr)$ is represented by
  $\lambda\cdot[C_+]$ or $\lambda\cdot[C_-]$, and
  \[
  f_3\bigl(\qform{1,-\lambda}h\bigr) = 2(\lambda\cdot[C_+]) =
  2(\lambda\cdot[C_-]).
  \]
  By \cite[(9.12)]{KMRT} we have
  \[
  2[C_+]=2[C_-] = 
  \begin{cases}
    0&\text{if $n\equiv0\bmod4$,}\\
    [D]&\text{if $n\equiv2\bmod4$.}
  \end{cases}
  \]
  The formula for $f_3\bigl(\qform{1,-\lambda}h\bigr)$ follows.

  Assume now $h$ has nontrivial discriminant. 
  The proof in this case follows the same pattern.
  Let $K/F$ be the
  discriminant field extension.   We consider the group scheme
  $R_{K/F}(\mathbb{G}_m)$, which is the Weil transfer of the
  multiplicative group. For every field $L$ containing $F$, the map
  \[
  \mu\in R_{K/F}(\mathbb{G}_m)(L)=(L\otimes_FK)^\times \mapsto
  e_3(\qform{1,-N_{L\otimes K/L}(\mu)} h_L)\in M^3_D(L)
  \]
  defines a degree~3 invariant of the quasi-trivial torus
  $R_{K/F}(\mathbb{G}_m)$ with values in the cycle module $M_D$.  By
  \cite[Th.~1.1]{MPT}, there is an element $u\in M_D^2(K)$ such that
  for any 
  field $L$ containing $F$ and any $\mu\in (L\otimes K)^\times$,
  \begin{equation*}
    \label{eq:ehypad2}
    e_3(\qform{1,-N_{L\otimes K/L}(\mu)} h_L) = \cores_{L\otimes
      K/L}(\mu\cdot u_{L\otimes K})\quad\text{in $M_D^3(L)$.}
  \end{equation*}
  It remains to show that $u=e_2(h_K)$. To prove this, we
  consider the field $L=K(t)$, where $t$ is an indeterminate. Since
  $e_1(h_{K(t)})=0$, the previous case applies. We thus get for any
  $\mu\in(K(t)\otimes_FK)^\times$
  \begin{equation*}
    \label{eq:ehypad3}
    N_{K(t)\otimes K/K(t)}(\mu)\cdot e_2(h_{K(t)}) =
    \cores_{K(t)\otimes K/K(t)}(\mu\cdot u_{K(t)\otimes K})
    \quad\text{in $M_D^3(K(t))$.}
  \end{equation*}
  Let $\iota$ be the nontrivial $F$-automorphism of $K$.
  The $K(t)$-algebra isomorphism $K(t)\otimes_FK\simeq K(t)\times K(t)$
  mapping $\alpha\otimes\beta$ to $(\alpha\beta,\alpha\iota(\beta))$
  yields an isomorphism $M_D^2(K(t)\otimes K)\simeq M_D^2(K(t))\times
  M_D^2(K(t))$ that carries $u_{K(t)\otimes K}$ to
  $(u_{K(t)},\iota(u)_{K(t)})$. Thus, for every $(\mu_1,\mu_2)\in
  K(t)^\times\times K(t)^\times$,
  \[
  \mu_1\mu_2\cdot e_2(h_{K(t)}) = \mu_1\cdot u_{K(t)} +\mu_2\cdot
  \iota(u)_{K(t)} \quad\text{in $M_D^3(K(t))$.}
  \]
  In particular, if $\mu_1=t$ and $\mu_2=1$ we get $t\cdot
  e_2(h_{K(t)})= t\cdot u_{K(t)}$, hence taking the residue for the
  $t$-adic valuation yields $e_2(h_K)=u$, proving the formula for
  $e_3\bigl(\qform{1,-N_{K/F}(\mu)}h\bigr)$. 

  To complete the proof, we compute
  $f_3\bigl(\qform{1,-N_{K/F}(\mu)}h\bigr)$. Let $C$ be the Clifford
  algebra of $\Ad_h$, so $[C]$ represents $e_2(h_K)$ and
  \[
  f_3\bigl(\qform{1,-N_{K/F}(\mu)}h\bigr)=2\cores_{K/F}(\mu\cdot[C]).
  \]
  By \cite[(9.12)]{KMRT} we have
  \[
  2[C]=
  \begin{cases}
    0&\text{if $n\equiv0\bmod4$,}\\
    [D_K]&\text{if $n\equiv2\bmod4$.}
  \end{cases}
  \]
  The formula for $f_3\bigl(\qform{1,-N_{K/F}(\mu)}h\bigr)$ follows by
  the projection formula.
\end{proof}

\subsection{Hermitian forms with an additive decomposition}

We now present another approach for computing the $f_3$-invariant, which does not rely on the computation of the Arason invariant. This leads to an explicit formula in a more general situation, which will be used in the proof of Theorem~\ref{f3.thm}:
\begin{prop}
  \label{prop:f3sum}
  Let $(V_1,h_1)$, \ldots, $(V_m,h_m)$ be hermitian modules of
  orthogonal type and even absolute rank $n_1$, \ldots, $n_m$ over
  $(D,\theta)$, and let $\lambda_1$, \ldots, $\lambda_m\in
  F^\times$. Let also $h=\qform{1,-\lambda_1}h_1\perp\ldots
  \perp\qform{1,-\lambda_m}h_m$. If $\sum_{i=1}^m\lambda_i\cdot
  e_1(h_i)=0$, then $h$ has trivial Clifford invariant, and 
  \[
  f_3(h) = \lambda_1^{n_1/2}\ldots \lambda_m^{n_m/2}\cdot[D].
  \]
\end{prop}

To prove this proposition, we need some preliminary results.
Let $(V,h_0)$ be a
hyperbolic module of orthogonal type over $(D,\theta)$. Recall from
\cite[(13.31)]{KMRT} the canonical map (``vector representation'')
\[
\chi\colon\Spin(h_0)\to \so(h_0).
\]
Since proper isometries have reduced norm~$1$, we also have the
inclusion
\[
i\colon \so(h_0) \to \SL(V).
\]

\begin{lem}
  \label{lem:Dynkindex}
  The following diagram, where $R$ is the Rost invariant, is
  commutative:
  \[
  \xymatrix{
  H^1(F,\Spin(h_0)) \ar[r]^{(i\circ\chi)_*} \ar[d]_{R} & H^1(F,\SL(V))
  \ar[d]^{R}\\
  H^3(F) \ar[r]^{2}& H^3(F)
  }
  \]
\end{lem}

\begin{proof}
  This lemma is just a restatement of the property that the Rost
  multiplier of the map $i\circ\chi$ is 2, see \cite[Ex.~7.15,
  p.~124]{GMS}. 
\end{proof}

We next recall from \cite[(29.27)]{KMRT} (see also \cite{GTW}) the
canonical description of the pointed set $H^1(F,\so(h_0))$. Define a
functor $\SSym(h_0)$ from $\Fields_F$ to the category of pointed sets
as follows: for any field $L$ containing $F$, set
\[
\SSym(h_0)(L) = \{(s,\lambda)\in \GL(V_L)\times L^\times \mid
\ad_{h_0}(s)=s \text{ and }\Nrd(s)=\lambda^2\},
\]
where the distinguished element is $(1,1)$. Let $F_s$ be a separable
closure of $F$ and let $\Gamma=\gal(F_s/F)$ be the Galois group. We
may identify $\SSym(h_0)(F_s)$ with the quotient
$\GL(V_{F_s})/\so((h_0)_{F_s})$ by mapping a class $a\cdot\so((h_0)_{F_s})$ to
$(a\ad_{h_0}(a),\Nrd(a))$ for $a\in \GL(V_{F_s})$. Therefore, we have
an exact sequence of pointed $\Gamma$-sets
\[
1\to \so((h_0)_{F_s}) \to \GL(V_{F_s}) \to \SSym(h_0)(F_s)\to 1.
\]
Since $H^1(F,\GL(V_{F_s}))=1$ by Hilbert's Theorem~90, the induced
exact sequence in Galois cohomology yields a canonical bijection
between $H^1(F,\so(h_0))$ and the orbit set of $\GL(V)$ on
$\SSym(h_0)(F)$. Abusing notation, we write simply $\SSym(h_0)$ for
$\SSym(h_0)(F)$. The orbits of $\GL(V)$ on $\SSym(h_0)$ are the
equivalence classes under the following relation:
\[
(s,\lambda)\sim(s',\lambda')\quad\text{if $s'=as\ad_{h_0}(a)$ and
  $\lambda'=\lambda\Nrd(a)$ for some $a\in\GL(V)$.}
\]
Therefore, we may identify
\[
H^1(F,\so(h_0)) = \SSym(h_0)/{\sim}.
\]

\begin{lem}
  \label{lem:RostSL}
  The composition $H^1(F,\so(h_0))\xrightarrow{i_*} H^1(F,\SL(V))
  \xrightarrow{R} H^3(F)$ maps the equivalence class of $(s,\lambda)$
  to $\lambda\cdot[A]$.
\end{lem}

\begin{proof}
  Let $\pi\colon\SSym(h_0)(F_s)\to F_s^\times$ be the projection
  $(s,\lambda)\mapsto\lambda$. We have a commutative diagram of
  pointed $\Gamma$-sets with exact rows:
  \[
  \xymatrix{1\ar[r]&\so((h_0)_{F_s}) \ar[r] \ar[d]_{i} & \GL(V_{F_s})
    \ar[r] \ar@{=}[d] & \SSym(h_0)(F_s) \ar[r] \ar[d]^{\pi}& 1\\
  1\ar[r]&\SL(V_{F_s})\ar[r]&\GL(V_{F_s})\ar[r]^{\Nrd}&
  {F_s^\times}\ar[r]& 1}
  \]
  This diagram yields the following commutative square in cohomology:
  \[
  \xymatrix{
  \SSym(h_0) \ar[r]\ar[d]_{\pi} & H^1(F,\so(h_0))\ar[d]^{i_*}\\
  F^\times\ar[r]& H^1(F,\SL(V))
  }
  \]
  On the other hand, the Rost invariant and the map $F^\times\to
  H^3(F)$ carrying $\lambda$ to $\lambda\cdot[A]$ fit in the following
  commutative diagram (see \cite[p.~437]{KMRT}):
  \[
  \xymatrix{
  F^\times\ar[rr]\ar[dr]&& H^1(F,\SL(V))\ar[dl]^{R}\\
  &H^3(F)&
  }
  \]
  The lemma follows.
\end{proof}

For the next statement, let $\partial\colon H^1(F,\so(h_0)) \to
{}_2H^2(F)$ be the connecting map in the cohomology exact sequence
associated to
\[
1\to\mu_2\to\Spin(h_0)\xrightarrow{\chi}\so(h_0)\to1.
\]
For any hermitian form $h$ of orthogonal type on $V$, there exists a
unique linear transformation $s\in\GL(V)$ such that
$h(x,y)=h_0(x,s^{-1}(y))$ for all $x$, $y\in V$, hence
$\ad_h=\Int(s)\circ \ad_{h_0}$ and $\ad_{h_0}(s)=s$. If the
discriminant of $h$ is trivial we have $\Nrd(s)\in F^{\times2}$, hence
there exists $\lambda\in F^\times$ such that $\lambda^2=\Nrd(s)$, and
we may consider $(s,\lambda)$ and $(s,-\lambda)\in \SSym(h_0)$. By
the main theorem of~\cite{GTW}, $\partial(s,\lambda)$ and
$\partial(s,-\lambda)$ are the Brauer 
classes of the two components of the Clifford algebra of
$\Ad_{h_0\perp-h}$, so if the Clifford invariant of $h$ is trivial we
have
\[
\{\partial(s,\lambda), \partial(s,-\lambda)\} = \{0,[D]\}.
\]

\begin{lem}
  \label{lem:altf3}
  With the notation above, we have $f_3(h)=\lambda\cdot[D]$ if
  $\partial(s,\lambda)=0$. 
\end{lem}

\begin{proof}
  By definition of $s$, the torsor in $H^1(F,\ort(h_0))$ corresponding
  to $h$ lifts to $(s,\lambda)\in H^1(F,\so(h_0))$. If
  $\partial(s,\lambda)=0$, then $(s,\lambda)$ lifts to some $\xi\in
  H^1(F,\Spin(h_0))$, and by definition of the invariants $e_3$ and
  $f_3$ we have
  \[
  e_3(h)=R(\xi)+F^\times\cdot[D]\in M^3_D(F) \qquad\text{and}\qquad
  f_3(h)=2R(\xi)\in H^3(F).
  \]
  Lemma~\ref{lem:Dynkindex} then yields
  $f_3(h)=R\circ(i\circ\chi)_*(\xi) = R\circ i_*(s,\lambda)$, and by
  Lemma~\ref{lem:RostSL} we have $R\circ
  i_*(s,\lambda)=\lambda\cdot[D]$. 
\end{proof}

In order to check the condition $\partial(s,\lambda)=0$ in
Lemma~\ref{lem:altf3}, the following observation is useful:
Suppose $(V_1,h_1)$ and $(V_2,h_2)$ are hermitian modules of
orthogonal type over $(D,\theta)$. The inclusions $V_i\hookrightarrow
V_1\perp V_2$ yield an $F$-algebra homomorphism $C(\Ad_{h_1})
\otimes_F C(\Ad_{h_2})\to C(\Ad_{h_1\perp h_2})$, which induces a
group homomorphism $\Spin(h_1)\times \Spin(h_2) \to \Spin(h_1\perp
h_2)$. This homomorphism fits into the following commutative diagram
with exact rows
\[
\xymatrix{
1\ar[r]&\mu_2\times\mu_2\ar[r]\ar[d]_{\prod} & \Spin(h_1)\times
\Spin(h_2) \ar[r]^{\chi_1\times\chi_2}\ar[d] & \so(h_1)\times
\so(h_2) \ar[r]\ar[d]^{\oplus}&1\\
1\ar[r]& \mu_2\ar[r]&\Spin(h_1\perp h_2) \ar[r]^{\chi}& \so(h_1\perp
h_2) \ar[r]&1
}
\]
The left vertical map is the product, and the right vertical map
carries $(g_1,g_2)$ to $g_1\oplus g_2$. The induced diagram in
cohomology yields the commutative square
\[
\xymatrix{
H^1(F,\so(h_1))\times H^1(F,\so(h_2))
\ar[r]^(.6){\partial_1\times \partial_2} \ar[d]_{\oplus} &
{}_2H^2(F)\times {}_2H^2(F)\ar[d]^{\prod}\\
H^1(F,\so(h_1\perp h_2)) \ar[r]^{\partial} & {}_2H^2(F)
}
\]
The following additivity property of the connecting maps $\partial$
follows: for $(s_1,\lambda_1)\in H^1(F,\so(h_1))$ and
$(s_2,\lambda_2)\in H^1(F,\so(h_2))$,
\begin{equation}
  \label{eq:additivity}
  \partial_1(s_1,\lambda_1) + \partial_2(s_2,\lambda_2) =
  \partial(s_1\oplus s_2,\lambda_1\lambda_2).
\end{equation}

\begin{proof}[Proof of Proposition~\ref{prop:f3sum}]
  Let $h_0=\qform{1,-1}h_1\perp\ldots\perp\qform{1,-1}h_m$, which is a
  hyperbolic form, and let $V=V_1^{\oplus2}\oplus\cdots\oplus
  V_m^{\oplus2}$ be the underlying vector space of $h$ and $h_0$. The
  linear transformation $s\in\GL(V)$ such that $h(x,y)=h_0(x,s^{-1}(y))$
  for all $x$, $y\in V$ is
  \[
  s=1\oplus\lambda_1^{-1}\oplus1\oplus\lambda_2^{-1}\oplus\cdots
  \oplus1\oplus\lambda_m^{-1}. 
  \]
  By the additivity property~\eqref{eq:additivity}, the connecting map
  \[\partial\colon H^1(F,\so(h_0))\to {}_2H^2(F)\] satisfies
  \[
  \partial(s,\lambda_1^{-n_1/2}\ldots \lambda_m^{-n_m/2})
  = \partial_1(\lambda_1^{-1},\lambda_1^{-n_1/2})
  +\cdots+\partial_m(\lambda_m^{-1}, \lambda_m^{-n_m/2}).
  \]
  A theorem of Bartels \cite[p.~283]{Bartels} (see also \cite{GTW})
  yields 
  $\partial_i(\lambda_i^{-1},\lambda_i^{-n_i/2}) = \lambda_i^{-1}\cdot
  e_1(h_i)$ for all $i$. Therefore, if $\sum_{i=1}^m\lambda_i\cdot
  e_1(h_i)=0$ 
  we have $f_3(h) = \lambda_1^{n_1/2}\ldots\lambda_m^{n_m/2}\cdot [D]$
  by Lemma~\ref{lem:altf3}.
\end{proof}

\subsection{Relative Arason invariant of hermitian forms of orthogonal
type} 
\label{subsec:relAra}

By using the Rost invariant, one may also define a relative
Arason invariant, in a broader context:

\begin{defi} Let $h_1$ and $h_2$ be two hermitian forms of orthogonal
  type over 
  $(D,\theta)$ such that their difference $h_1+(-h_2)$ has even
  relative rank, trivial discriminant, and trivial Clifford
  invariant. Their relative Rost invariant
  is defined by
  \[
  e_3(h_1/h_2)=e_3\bigr(h_1\perp(-h_2)\bigl)\in M_D^3(F).
  \]
\end{defi}

In particular, if both $h_1$ and $h_2$ have even relative rank,
trivial discriminant, and trivial Clifford invariant, then
$e_3(h_1/h_2)=e_3(h_1)+e_3(h_2)=e_3(h_1)-e_3(h_2)$.

\begin{rem} 
  Under the conditions of this definition, one may check
  that the involution $\ad_{h_2}$ corresponds to a torsor which can be
  lifted to a $\Spin(\Ad_{h_1})$ torsor (see~\cite[\S 3.5]{Hyd}). As
  explained in~\cite[Lemma 3.6]{BFP98}, the relative Arason invariant
  $e_3(h_1/h_2)$ coincides with the class in $M_D^3(F)$ of the image
  of this torsor under the Rost invariant of $\Spin(\Ad_{h_1})$.
\end{rem}

Combining the properties of the Arason invariant recalled in
Lemma~\ref{Arason-herm.lem} and the computation of
Proposition~\ref{compute.prop}, we obtain:

\begin{cor}
  \label{compute.cor}
  \begin{enumerate}
  \item[(i)] Let $h$ be a hermitian form of orthogonal type with even
    absolute rank, and let 
    $K/F$ be the discriminant quadratic extension.  For any $\mu\in
    K^\times$, the relative Arason invariant $e_3(\qform{N_{K/F}(\mu)}
    h/h)$ is well-defined, and
    \[
    e_3(\qform{N_{K/F}(\mu)} h/h)=\cores_{K/F}(\mu\cdot e_2(h_K)).
    \]
  \item[(ii)] Let $h_1$ and $h_2$ be two hermitian forms of orthogonal type
    with even absolute 
    rank and trivial discriminant. We have
    \[
    e_3(h_1\perp\qform{\lambda}h_2/h_1\perp h_2)=
    e_3(\qform{\lambda}h_2/h_2)=\lambda\cdot e_2(h_2).
    \]
  \end{enumerate}
\end{cor}

\subsection{Arason and $f_3$ invariants of orthogonal involutions}
\label{subsec:Arasoninvol}

Let $(A,\sigma)$ be an algebra with orthogonal involution,
Brauer-equivalent to the division algebra $D$ over $F$. We pick an
involution $\theta$ on $D$, so that $(A,\sigma)$ can be represented as
the adjoint $(A,\sigma)\simeq \Ad_h$ of some hermitian module
$(V,h)$ over $(D,\theta)$. The co-index of $A$, which is the dimension over $D$ of the 
module $V$, is equal to the relative rank of $h$.  If the form $h$ has even relative rank, trivial
discriminant, and trivial Clifford invariant, then its Arason
invariant is well-defined. Moreover, by Lemma~\ref{Arason-herm.lem},
we have $e_3(h)=e_3(\lambda h)$ for any $\lambda\in F^\times$, and, as
explained in~\cite[Prop 3.8]{BFP98}, $e_3(h)$ does not depend on the
choice of $\theta$. Therefore, we get a well-defined Arason
invariant for the
involution $\sigma$, provided the algebra $A$ has even co-index,
i.e. $\deg(A)/\ind(A)=\deg(A)/\deg(D)\in 2\mz$, and the involution
$\sigma$ has trivial discriminant and trivial Clifford invariant:
$$
e_3(\sigma)=e_3(h)\in M^3_A(F)=M^3_D(F).
$$

\begin{rems}
  \begin{enumerate}
  \item Under the assumptions above on $(A,\sigma)$, one may also check
    that the algebra $A$ carries a hyperbolic orthogonal
    involution $\sigma_0$, and the Arason invariant $e_3(\sigma)$ can
    be defined directly in terms of the Rost invariant of the group
    $\Spin(A,\sigma_0)$, see~\cite[\S 3.5]{Hyd}.
  \item Similarly, we may also define a relative Arason invariant
    $e_3(\sigma_1/\sigma_2)$ if the involutions $\sigma_1$ and
    $\sigma_2$ both have trivial discriminant and trivial Clifford
    invariant. But we cannot relax those assumptions, as we did for
    hermitian forms. Indeed, if $e_2(h_2)=e_2(\ad_{h_2})$ is not
    trivial, then $e_3(\qform{\lambda}h_1/h_2)$ and $e_3(h_1/h_2)$ are
    generally different, as Corollary~\ref{compute.cor} shows.
  \end{enumerate}
\end{rems}

In the setting above, we may also define an $f_3$-invariant by
$f_3(\sigma)=f_3(h)$, or equivalently:

\begin{defi}
\label{f3.inv}
  Let $(A,\sigma)$ be an algebra with orthogonal involution. We assume
  $A$ has even co-index, and $\sigma$ has trivial discriminant and
  trivial Clifford invariant.  We define $f_3(\sigma)\in
  F^\times\cdot[A]\subset {}_2H^3(F)$ by 
  $f_3(\sigma)=2c$, where $c$ is any representative of the Arason
  invariant $e_3(\sigma)\in M^3_A(F)$. 
\end{defi}

\begin{rem}
(i) If $A$ is split, then $F^\times \cdot [A]=\{0\}$, and
  $f_3(\ad_\varphi)=0$  for all quadratic forms $\varphi\in
  I^3(F)$. This also follows from the fact that $e_3(\varphi)\in
  {}_2H^3(F)$.

 (ii) Using the same process, one may define an invariant $f_3^{16}$ from Garibaldi's invariant $e_3^{16}$, and from Bermudez-Ruozzi's generalization (see~\cite{Ga},~\cite{BR}). This invariant has values in ${}_2H^3(F)$, but need not have values in $F^\times\cdot [A]$ in general. 
  \end{rem} 

\begin{ex} 
  Let $Q$ be a quaternion algebra over $F$, and consider the algebra with involution $(A,\sigma)=(Q,\rho)\otimes\Ad_\varphi$, where $\rho$
  is an orthogonal involution with discriminant $\delta \cdot
  F^{\times 2}\in F^\times/F^{\times 2}$, and $\varphi$ is a
  even-dimensional quadratic form with trivial discriminant. We have
  $e_3(\sigma)=\delta\cdot e_2(\varphi)\mod F^\times \cdot [Q]$, and
  $f_3(\sigma)=0$. Indeed, since the restriction map $M_Q^3(F)\ra
  M^3_Q(F_Q)=H^3(F_Q)$ is injective, it is enough to check the formula
  in the split case, where it follows from a direct computation.
\end{ex}

The computation in Proposition~\ref{compute.prop} can be again
rephrased as follows:

\begin{cor}
  \label{compute-inv.cor}
  Let $(A,\sigma)$ be a central simple $F$-algebra of even degree~$n$
  with orthogonal 
  involution, and let $K/F$ be the discriminant quadratic
  extension. For any $\mu\in K^\times$, the algebra with involution
  $\Ad_\qform{1,-N_{K/F}(\mu)}\otimes (A,\sigma)$ has even co-index,
  trivial discriminant and trivial Clifford invariant. Its
  Arason invariant is given by
  \[
  e_3(\ad_\qform{1,-N_{K/F}(\mu)}\otimes
  \sigma)=\cores_{K/F}\bigl(\mu\cdot e_2(\sigma_K)\bigr),
  \]
  and 
  \[
  f_3(\ad_\qform{1,-N_{K/F}(\mu)}\otimes
  \sigma)=
  \begin{cases}
    0&\text{if $n\equiv 0\bmod4$,}\\
    N_{K/F}(\mu)\cdot[A]&\text{if $n\equiv2\bmod4$.}
  \end{cases}
  \]
  In particular, if $\sigma$ has trivial discriminant, we have for any
  $\lambda\in F^\times$
  \[
  e_3(\ad_\qform{1,-\lambda}\otimes \sigma)=\lambda\cdot
  e_2(\sigma)
  \]
  and
  \[
  f_3(\ad_\qform{1,-\lambda}\otimes \sigma)=
  \begin{cases}
    0&\text{if $n\equiv0\bmod4$,}\\
    \lambda\cdot[A]&\text{if $n\equiv2\bmod4$.}
  \end{cases}
  \]
\end{cor}

Hence, the formula given in~\cite[Th.~5.5]{QT:Arason} for algebras of
degree $8$ is actually valid in arbitrary degree.

With this in hand, one may easily check that the $f_3$ invariant is
trivial up to degree $10$. Indeed, since the co-index of the algebra is
supposed to be even, the algebra is possibly non-split only when its
degree is divisible by~$4$. In degree $4$, any involution with trivial
discriminant 
and Clifford invariant is hyperbolic, hence has trivial invariants. In
degree $8$, any involution with trivial
discriminant 
and Clifford invariant admits a decomposition as in 
Corollary~\ref{compute-inv.cor} by~\cite[Th.~5.5]{QT:Arason}, hence
its $f_3$ invariant is trivial. In degree $12$, one may construct
explicit examples of $(A,\sigma)$ with $f_3(\sigma)\not=0$ as follows. Suppose
$E$ is a central simple $F$-algebra of degree~$4$. Recall from
\cite[\S10.B]{KMRT} that the second $\lambda$-power $\lambda^2E$ is a
central simple $F$-algebra of degree~$6$, which carries a canonical
involution $\gamma$ of orthogonal type with trivial discriminant, and
which is Brauer-equivalent to $E\otimes_FE$. 
\begin{cor}
  \label{f3.cor}
  Let $E$ be a central simple $F$-algebra of degree and
  exponent~$4$. Pick an indeterminate $t$, and consider the algebra 
  with involution \[(A,\sigma)=\Ad_{\qform{1,-t}}\otimes
  (\lambda^2E,\gamma)_{F(t)}.\] We
  have 
  \[
  f_3(\sigma)=t\cdot[A]\neq0\in H^3(F(t)).
  \]
\end{cor}

\begin{proof}
  The formula $f_3(A,\sigma)=t\cdot[A]$ readily follows from
  Corollary~\ref{compute-inv.cor}. The algebra $E$ has exponent $4$,
  therefore  
  $[A]=[E\otimes_F E]\not=0$. 
  Since $t$ is an indeterminate, we get $t\cdot[A]\neq0$.
\end{proof}
 
\section{Additive decompositions in degree $12$}
\label{sec:adddec12}

In the next three sections, we concentrate on degree $12$ algebras
$(A,\sigma)$ with orthogonal involution of trivial discriminant and
trivial Clifford invariant.  
The main result of this section is Theorem~\ref{thm:adddec12}, which
generalizes a theorem of Pfister on $12$-dimensional quadratic forms. 

\subsection{Additive decompositions} 

Given three algebras with involution, $(A,\sigma)$, $(A_1,\sigma_1)$ and
$(A_2,\sigma_2)$, we say that $(A,\sigma)$ is a \emph{direct sum} of
$(A_1,\sigma_1)$ and $(A_2,\sigma_2)$, and we write
$$
(A,\sigma)\in (A_1,\sigma_1)\orthsum(A_2,\sigma_2),
$$ 
if there exist a division algebra with involution $(D,\theta)$ and
hermitian modules $(V_1,h_1)$ and $(V_2,h_2)$ over $(D,\theta)$, which
are both hermitian or both skew-hermitian, such that
$(A_1,\sigma_1)=\Ad_{h_1}$, $(A_2,\sigma_2)=\Ad_{h_2}$ and
$(A,\sigma)=\Ad_{h_1\perp h_2}$. In particular, this implies $A$,
$A_1$ and $A_2$ are all three Brauer-equivalent to $D$, and the
involutions $\sigma$, $\sigma_1$ and $\sigma_2$ are of the same type.  
This notion of direct sum for algebras with involution was introduced
by Dejaiffe in~\cite{Dejaiffe}. As explained there, the algebra with
involution $(A,\sigma)$ is generally not uniquely determined by the
data of the two summands $(A_1,\sigma_1)$ and $(A_2,\sigma_2)$.  Indeed,
multiplying the hermitian forms $h_1$ and $h_2$ by a scalar does not
change the adjoint involutions, so the adjoint of ${\lambda_1
  h_1\perp\lambda_2 h_2}$ also is a direct sum of $(A_1,\sigma_1)$ and
$(A_2,\sigma_2)$ for any $\lambda_1$, $\lambda_2\in F^\times$. 
If one of the two summands, say $(A_1,\sigma_1)=(A_1,\hyp)$ is hyperbolic, then all hermitian forms similar to $h_1$ actually are isomorphic to $h_1$. Hence in this case, there is a unique direct sum, and we may write 
\[(A,\sigma)=(A_1,\hyp)\boxplus(A_2,\sigma_2).\]

The cohomological invariants we consider, when defined, have the
following additivity property: 

\begin{prop}
  \label{prop:additinv}
  Suppose $\sigma$, $\sigma_1$, $\sigma_2$ are orthogonal involutions
  such that $(A,\sigma)\in(A_1,\sigma_1)\boxplus(A_2,\sigma_2)$. We have:
  \begin{enumerate}
  \item[(i)]
  $\deg A=\deg A_1+\deg A_2$.
  \item[(ii)]
  If $\deg A_1\equiv\deg A_2\equiv0\bmod2$, then
  $e_1(\sigma)=e_1(\sigma_1)+e_1(\sigma_2)$.
  \item[(iii)]
  If $\deg A_1\equiv\deg A_2\equiv0\bmod2$ and
  $e_1(\sigma_1)=e_1(\sigma_2)=0$, then
  \[e_2(\sigma)=e_2(\sigma_1)+e_2(\sigma_2).\]
  \item[(iv)]
  If the co-indices of $A_1$ and $A_2$ are even and
  $e_i(\sigma_1)=e_i(\sigma_2)=0$ for $i=1$, $2$, then
  \[e_3(\sigma)=e_3(\sigma_1) + e_3(\sigma_2)\mbox{ and }f_3(\sigma)=
  f_3(\sigma_1)+f_3(\sigma_2).\]
  \end{enumerate}
\end{prop}

\begin{proof}
  Assertion (i) is clear by definition, and (ii) was established by
  Dejaiffe~\cite[Prop.~2.3]{Dejaiffe}. Assertion~(iii) follows from
  \cite[\S\,3.3]{Dejaiffe} (see also the proof of the ``Orthogonal Sum
  Lemma'' in \cite[\S3]{Garihyp}).
  To prove~(iv), let $D$ be the division algebra Brauer-equivalent to
  $A$, $A_1$, and $A_2$, and let $\theta$ be an $F$-linear involution
  on $D$. We may find hermitian forms of orthogonal type $h_1$, $h_2$
  over $(D,\theta)$ such that $(A_i,\sigma_i)\simeq\Ad_{h_i}$ for
  $i=1$, $2$, and $(A,\sigma)\simeq\Ad_{h_1\perp h_2}$. By
  Lemma~\ref{Arason-herm.lem}(ii) we have \[e_3(h_1\perp h_2)=e_3(h_1)
  + e_3(h_2).\]
  By definition of the $e_3$-invariant of orthogonal
  involutions (see \S\,\ref{subsec:Arasoninvol}), $e_3(\sigma)$
  (resp.\ $e_3(\sigma_i)$ for $i=1$, $2$) is represented by
  $e_3(h_1\perp h_2)$ (resp.\ $e_3(h_i)$), hence
  $e_3(\sigma) = e_3(\sigma_1)+e_3(\sigma_2)$. Likewise, the
  additivity of $e_3$ induces $f_3(\sigma) =
  f_3(\sigma_1)+f_3(\sigma_2)$, by definition of the $f_3$-invariant
  (see~\ref{f3.inv}).
\end{proof}

By a result of Pfister~\cite[p.123-124]{Pfister}, any $12$-dimensional quadratic form $\varphi$ in
$I^3F$ decomposes as $\varphi=\qform{\alpha_1}n_1\perp
\qform{\alpha_2}n_2\perp\qform{ \alpha_3}n_3$, where $n_i$ is a
$2$-fold Pfister form and $\alpha_i\in F^\times$, for $1\leq i\leq
3$. This can be rephrased as 
$$\Ad_\varphi\in\Ad_{n_1}\boxplus\Ad_{n_2}\boxplus\Ad_{n_3},$$
where each summand $\Ad_{n_i}$ has degree  $4$ and discriminant $1$. 
We now extend this result to the non-split case.

\begin{thm}
  \label{thm:adddec12}
  Let $(A,\sigma)$ be a central simple $F$-algebra of degree~$12$ with
  orthogonal involution. Assume the discriminant and the Clifford
  invariant of $\sigma$ are trivial. There is a central simple
  $F$-algebra $A_0$ of degree~$4$ and orthogonal involutions
  $\sigma_1$, $\sigma_2$, $\sigma_3$ of trivial discriminant on $A_0$
  such that
  \[
  (A,\sigma)\in (A_0,\sigma_1)\orthsum (A_0,\sigma_2)
  \orthsum(A_0,\sigma_3).
  \]
\end{thm}

Note that since $\deg A_0=4$ we have $e_2(\sigma_i)=0$ if and only if
$\sigma_i$ is hyperbolic (see~\cite[Th.~3.10]{Hyd}); therefore, even
when the index of $A$ is~$2$ we \emph{cannot} use
Proposition~\ref{prop:additinv}(iv) to compute $e_3(\sigma)$ (unless
each $\sigma_i$ is hyperbolic).

\begin{proof}[Proof of Theorem~\ref{thm:adddec12}]
  The index of $A$ is a power of~$2$ since $2[A]=0$ in $\Br(F)$, and
  it divides~$\deg A=12$, so $\ind A=1$, $2$ or $4$. As we just
  pointed out, the  
  index~$1$ case is Pfister's theorem. We consider separately the two
  remaining cases. 
  
    If $\ind A=2$, we have $(A,\sigma)=\Ad_h$ for some skew-hermitian
  form $h$ of relative rank~$6$ over a quaternion division algebra
  $(Q,\ba)$ with its canonical involution. Let $q_1\in Q$ be a
  nonzero pure
  quaternion represented by $h$, and write $h=\qform{q_1}\perp
  h'$. Over the quadratic extension $K_1=F(q_1)$, the algebra $Q$
  splits and the form $\qform{q_1}$ becomes hyperbolic (because its
  discriminant becomes a square). Therefore, $h_{K_1}$ and
  $h'_{K_1}$ are Witt-equivalent, and $(\ad_{h'})_{K_1}$ is adjoint to
  a $10$-dimensional form $\varphi$. The discriminant and Clifford
  invariant of $\sigma$ are 
  trivial, hence $\varphi\in I^3K_1$. Since there is no anisotropic
  $10$-dimensional quadratic 
  forms in $I^3$ (see \cite[Th.~8.1.1]{Kahn}), it follows that
  $h'_{K_1}$ is isotropic, 
  hence  by
  \cite[Prop., p.~382]{QT}, $h'=\qform{-\lambda_1q_1}\perp k$ for some
  $\lambda_1\in 
  F^\times$ and some skew-hermitian form $k$ of relative rank~$4$. We thus have
  \[
  h=\qform{q_1}\qform{1,-\lambda_1}\perp k,
  \]
  and computation shows that
  $e_1\bigl(\qform{q_1}\qform{1,-\lambda_1}\bigr)=0$. Therefore,
  $e_1(k)=0$, and
  $e_2(k)=e_2\bigl(\qform{q_1}\qform{1,-\lambda_1}\bigr)$ because
  $e_2(h)=0$. Now, let $q_2\in Q$ be a nonzero pure quaternion
  represented by $k$, and let $K_2=F(q_2)$, so $k=\qform{q_2}\perp k'$
  for some skew-hermitian form $k'$ of relative rank~$3$. The forms
  $k_{K_2}$ and $k'_{K_2}$ are Witt-equivalent, and $(\ad_{k'})_{K_2}$
  is adjoint to a $6$-dimensional form $\psi\in I^2K_2$, i.e., to an
  Albert form $\psi$. We have
  \[
  e_2(\psi)=e_2(k')_{K_2} = e_2(k)_{K_2} =
  e_2\bigl(\qform{q_1}\qform{1,-\lambda_1}\bigr)_{K_2},
  \]
  hence the index of $e_2(\psi)$ is at most~$2$, and it follows that
  $\psi$ is isotropic. Therefore, $k'_{K_2}$ is isotropic, and
  $k'=\qform{-\lambda_2q_2}\perp\ell$ for some $\lambda_2\in
  F^\times$ and some skew-hermitian form $\ell$ of relative
  rank~$2$. Thus, we have
  \[
  h=\qform{q_1}\qform{1,-\lambda_1} \perp
  \qform{q_2}\qform{1,-\lambda_2} \perp\ell.
  \]
  Since $e_1\bigl(\qform{q_1}\qform{1,-\lambda_1}\bigr) =
  e_1\bigl(\qform{q_2}\qform{1,-\lambda_2}\bigr)=0$ and $e_1(h)=0$, we
  also have $e_1(\ell)=0$. We thus obtain the required decomposition,
  with
  \[
  A_0=M_2(Q),\qquad
  \sigma_1=\ad_{\qform{q_1}\qform{1,-\lambda_1}},\qquad
  \sigma_2=\ad_{\qform{q_2}\qform{1,-\lambda_2}},\qquad
  \sigma_3 = \ad_\ell.
  \]

  Suppose now $\ind A=4$, and let $D$ be the division algebra of
  degree~$4$ Brauer-equivalent to $A$. By~\cite[Th.~3.1]{GQ:deg12},
  there 
  exists a quadratic extension $K$ of $F$ such that $(A,\sigma)_K$ is
  hyperbolic. The co-index of $A_K$ is therefore even, so the index of
  $A_K$ is $2$, hence we may identify $K$ with a subfield of
  $D$. The following construction is inspired by the
  Parimala--Sridharan--Suresh exact sequence in Appendix~2 of
  \cite{BFP}. We have $D=\widetilde D\oplus D'$, where $\widetilde D$
  is the 
  centralizer of $K$ in $D$ and, writing $\iota$ for the nontrivial
  $F$-automorphism of $K$,
  \[
  D'=\{x\in D\mid xy=\iota(y)x\text{ for all $y\in K$}\}.
  \]
  Let $\theta$ be an orthogonal involution on $D$ that fixes $K$ (such
  involutions exist by \cite[(4.14)]{KMRT}). We may represent
  $(A,\sigma)=(\End_DV,\ad_h)$ for some hermitian form $h$ of relative
  rank~$3$ 
  over $(D,\theta)$. In view of the decomposition $D=\widetilde
  D\oplus D'$, we have for $x$, $y\in V$
  \[
  h(x,y)=\widetilde h(x,y) + h'(x,y) \qquad\text{with $\widetilde
    h(x,y)\in\widetilde D$ and $h'(x,y)\in D'$.}
  \]
  Since $h$ is a hermitian form over $(D,\theta)$, it follows that
  $\widetilde h$ is a hermitian form on $V$ viewed as a $\widetilde
  D$-vector space, with respect to the restriction of $\theta$ to
  $\widetilde D$. Clearly, $\End_DV\subset \End_{\widetilde D}V$. We
  may also embed $K$ into $\End_{\widetilde D}V$ by identifying
  $\alpha\in K$ with the scalar multiplication $x\mapsto x\alpha$ for
  $x\in V$. Thus, we have a $K$-algebra homomorphism
  \[
  (\End_{D}V)\otimes_FK\to \End_{\widetilde D}V.
  \]
  This homomorphism is injective because the left side is a simple
  algebra, hence it is an isomorphism by dimension count. For
  $f\in\End_DV$ we have $\ad_h(f)=\ad_{\widetilde h}(f)$, so the
  isomorphism preserves the involution, and therefore
  $(\Ad_h)_K=\Ad_{\widetilde h}$. Since $\sigma$ becomes hyperbolic
  over $K$, the form $\widetilde h$ is hyperbolic. Therefore, there is
  an $h$-orthogonal base of $V$ consisting of $\widetilde h$-isotropic
  vectors, which yields a diagonalization
  \[
  h=\qform{a_1,a_2,a_3} \qquad\text{with $a_1$, $a_2$, $a_3\in D'\cap
    \Sym(\theta)$.}
  \]
  We thus have $(A,\sigma)\in(D,\sigma_1)\orthsum(D,\sigma_2)
  \orthsum(D,\sigma_3)$ with $\sigma_i=\Int(a_i^{-1})\circ\theta$ for
  $i=1$, $2$, $3$. To complete the proof, we show that the
  discriminant of each $\sigma_i$ is trivial. Recall from
  \cite[(7.2)]{KMRT} that the discriminant is the square class of any
  skew-symmetric unit. Let $\alpha\in K^\times$
  be such that $\iota(\alpha)=-\alpha$. Since $a_i\in D'$ we have
  $\sigma_i(\alpha)=-\alpha$, so
  $\disc\sigma_i=\Nrd_D(\alpha)=N_{K/F}(\alpha)^2$.
\end{proof}

Recall that a central simple algebra of degree~$4$ with orthogonal
involution $(A_0,\sigma_0)$ of trivial discriminant decomposes as
$(A_0,\sigma_0)\simeq (Q,\ba)\otimes(H,\ba)$ where the quaternion
algebras $Q$, $H$ are the two components of the Clifford algebra
$C(A_0,\sigma_0)$ (see~\cite[(15.12)]{KMRT}). Therefore,
Theorem~\ref{thm:adddec12} can be rephrased as follows:

\begin{cor}
  \label{add-dec.cor}
  Let $(A,\sigma)$ be a central simple algebra of degree~$12$ with
  orthogonal involution of trivial discriminant and Clifford
  invariant. There exist quaternion $F$-algebras $Q_i$, $H_i$ for
  $i=1$, $2$, $3$ such that $[A]=[Q_i]+[H_i]$ for $i=1$, $2$, $3$,
  $[H_1]+[H_2]+[H_3]=0$, and
  \[
  (A,\sigma)\in \bigl((Q_1,\ba)\otimes(H_1,\ba)\bigr) \orthsum
  \bigl((Q_2,\ba)\otimes(H_2,\ba)\bigr) \orthsum
  \bigl((Q_3,\ba)\otimes(H_3,\ba)\bigr).
  \]
\end{cor}

\begin{proof}
  Theorem~\ref{thm:adddec12} yields orthogonal involutions $\sigma_1$,
  $\sigma_2$, $\sigma_3$ of trivial discriminant on the central simple
  $F$-algebra $A_0$ of degree~$4$ Brauer-equivalent to $A$ such that
  \[
  (A,\sigma)\in
  (A_0,\sigma_1)\orthsum(A_0,\sigma_2)\orthsum(A_0,\sigma_3).
  \]
  Each $(A_0,\sigma_i)$ has a decomposition
  \[
  (A_0,\sigma_i)\simeq (Q_i,\ba)\otimes(H_i,\ba)
  \]
  for some quaternion $F$-algebras $Q_i$, $H_i$ such that 
  \[
  e_2(\sigma_i) = [Q_i]+\{0,[A]\} = [H_i]+\{0,[A]\} \in M^2_A(F).
  \]
  From this decomposition, it follows that 
  \[
  [Q_i]+[H_i] = [A_0] = [A] \qquad\text{for $i=1$, $2$, $3$.}
  \]
  Moreover, we have $\sum_{i=1}^3e_2(\sigma_i)=e_2(\sigma)$
  by Proposition~\ref{prop:additinv}, and $e_2(\sigma)=0$, so
  \[
  \sum_{i=1}^3[Q_i]  = \bigl(\sum_{i=1}^3[H_i]\bigr) +[A]\in\{0,[A]\}
  \subset{}_2H^3(F).
  \]
  Therefore, interchanging $Q_i$ and $H_i$ if necessary, we may assume
  $\sum_{i=1}^3[H_i]=0$.
\end{proof}

\begin{ex} 
\label{split.ex}
For any quaternion algebra $Q$ with norm form $n_Q$, we have
$\Ad_{n_Q}\simeq (Q,\ba)\otimes (Q,\ba)$, see for instance
\cite[(11.1)]{KMRT}. Therefore, if $A$ is split, and $\sigma$ is
adjoint to the $12$-dimensional form
$\varphi=\qform{\alpha_1}n_1\perp\qform{\alpha_2}n_2\perp
\qform{\alpha_3}n_3$,  
where $n_i$ is the norm form of a quaternion algebra $Q_i$ for $1\leq
i\leq 3$, then  
$$(A,\sigma)\in\boxplus_{i=1}^3(Q_i,\ba)\otimes (Q_i,\ba).$$

Conversely, any decomposition of a split $(A,\sigma)=\Ad_\varphi$ as in the corollary corresponds to a decomposition 
$\varphi=\qform{\alpha_1}n_1\perp\qform{\alpha_2}n_2\perp
\qform{\alpha_3}n_3$, where $n_i$ is the norm form of 
$Q_i\simeq H_i$.   
\end{ex}

\begin{rem}
  \label{GaQ.rem}
  In~\cite{GQ:deg12}, it is proved that any $(A,\sigma)$ of degree
  $12$ with trivial discriminant and trivial Clifford invariant can be
  described as a quadratic extension of some degree $6$ central simple
  algebra with unitary involution $(B,\tau)$, with discriminant
  algebra 
  Brauer-equivalent to $A$. This algebra $(B,\tau)$ can be described
  from the above additive decomposition as follows.  
  Since 
  $\sum_{i=1}^3[H_i]=0$, the algebras $H_i$ have a common quadratic
  subfield $K$, see~\cite[Th.~III.4.13]{Lam}. All three products
  $(Q_i,\ba)\otimes (H_i,\ba)$ are 
  hyperbolic over $K$, so $\sigma_K$ is hyperbolic. Moreover, as
  observed in~\cite[Ex. 1.3]{GQ:deg12}, the tensor product
  $(Q_i,\ba)\ot(H_i,\ba)$ is a quadratic extension of
  $(Q_i,\ba)\otimes (K,\ba)$. Therefore, $(A,\sigma)$ is a quadratic
  extension of some $(B,\tau)\in \orthsum_{i=0}^3(Q_i,\ba)\otimes
  (K,\ba)$, and the discriminant algebra of $(B,\tau)$ is
  Brauer-equivalent to $[Q_1]+[Q_2]+[Q_3]=[Q]$.  
  Note that in the case
  where $\ind A=4$, we use in our proof the main result
  of~\cite{GQ:deg12}, which guarantees the existence of a quadratic
  extension $K$ such that $(A,\sigma)_K$ is hyperbolic.  
  But for $\ind A\neq 4$, our proof is independent, and does not use
  the existence of an open 
  orbit of a half-spin representation as in~\cite[p.~1220]{GQ:deg12}. 
\end{rem}

\subsection{Decomposition groups of $(A,\sigma)$}

Until the end of this section, $(A,\sigma)$ denotes a central simple
$F$-algebra of degree $12$ with an orthogonal involution of trivial
discriminant and trivial Clifford algebra.

\begin{defi}
  \label{defi:decsubgrp}
  Given an additive decomposition as in Corollary~\ref{add-dec.cor}
  \[
  (A,\sigma) \in \orthsum_{i=1}^3
  \bigl((Q_i,\ba)\otimes_F(H_i,\ba)\bigr)
  \qquad\text{with $\sum_{i=1}^3[H_i]=0$,}
  \]
  the subset
  \[
  U=\{0,[A],[Q_1], [H_1], [Q_2], [H_2], [Q_3],
  [H_3]\}\subset{}_2\Br(F)
  \]
  is called a \emph{decomposition group} of $(A,\sigma)$. It is
  indeed the subgroup of ${}_2\Br(F)$ generated by $[Q_1]$, $[Q_2]$,
  and $[Q_3]$, since $[A]=[Q_1]+[Q_2]+[Q_3]$ and $[H_i]=[A]+[Q_i]$ for
  $i=1$, $2$, $3$. 
\end{defi}

As the following examples show, a given algebra with
  involution $(A,\sigma)$ may admit several additive decompositions,
  corresponding to different decomposition groups, possibly not all of
  the same cardinality.

\begin{ex}
\label{ex:split}
Assume $A$ is split. Since $[A]=0$, we have $[H_i]=[Q_i]$ for all $i$. Hence all decomposition groups of $(A,\sigma)$ have order dividing $4$. 

Consider three quaternion division algebras $[Q_1]$, $[Q_2]$ and
$[Q_3]$ such that $[Q_1]+[Q_2]=[Q_3]$. By the ``common slot
  lemma'' \cite[Th.~III.4.13]{Lam}, there exist $a$, $b_1$, $b_2\in
F^\times$ such that $Q_i=(a,b_i)$ for $i=1,2$ and $Q_3=(a,b_1b_2)$. An
easy computation then shows that the norm forms of $Q_1$,
  $Q_2$, $Q_3$, respectively denoted by $n_1$, $n_2$, $n_3$, satisfy 
$n_1-n_2=\qform{b_2}n_3$ in the Witt group of $F$.  
Hence, extending scalars to a rational function fields in two variables over $F$, one may find scalars $\alpha_i$ for $1\leq i\leq 3$ such that the form 
$$\varphi=
\qform{\alpha_1}n_1\perp\qform{\alpha_2}n_2\perp\qform{\alpha_3}n_3$$
is either anisotropic, or isotropic and non-hyperbolic, or hyperbolic.  
By Example~\ref{split.ex}, in all three cases,
$\{0,[Q_1],[Q_2],[Q_3]\}$ is a decomposition group of order $4$ for
the involution $\sigma=\ad_\varphi$.  

On the other hand, the adjoint involution of an isotropic or
a hyperbolic form also has smaller decomposition groups, as we now
proceed to show.  
If the involution $\sigma$ is isotropic, it is adjoint to a quadratic
form $\varphi$ which is Witt-equivalent to a $3$-fold Pfister form
$\pi_3$. Let $Q$ be a quaternion algebra such that the norm form $n_Q$
is a subform of $\pi_3$. There exists $\alpha_1$, $\alpha_2\in
F^\times$ such that $\varphi=\qform{\alpha_1,\alpha_2}\otimes
n_Q\perp2\mh$, (where $\mh$ denotes the hyperbolic form)
hence $\{0,[Q]\}$ is a decomposition group of 
$\sigma=\ad\varphi$.  
If in addition $\sigma$, hence $\pi_3$, is hyperbolic, we may choose $[Q]=0$. 
\end{ex} 

\begin{ex}
\label{ex:ind2}
Assume now that $A=M_6(Q)$ has index $2$. Since $0\not =[Q]\in U$, all decomposition groups $U$ have order $2$, $4$ or $8$. 

If $\sigma$ is isotropic, then it is Witt-equivalent to a degree $8$
algebra with involution $(M_4(Q),\sigma_0)$ that has trivial
discriminant and trivial Clifford invariant, so
\[
(A,\sigma) = (M_4(Q),\sigma_0)\boxplus(M_2(Q),\hyp) =
(M_4(Q),\sigma_0)\boxplus \bigl((M_2(F),\ba)\otimes(Q,\ba)\bigr).
\]
Because $(M_4(Q),\sigma_0)$ has trivial discriminant and Clifford invariant,
by~\cite[Th.~5.2]{QT:Arason} we may find $\lambda$, $\mu\in F^\times$
and an orthogonal involution $\rho$ on $Q$ such that
\[
(M_4(Q),\sigma_0)\simeq\Ad_{\pform{\lambda,\mu}}\otimes (Q,\rho).
\]
Let $Q_1$ and $H_1$ be the two components of the Clifford algebra of
$\Ad_{\pform{\mu}}\otimes(Q,\rho)$. Then
\[
\Ad_{\pform{\mu}}\otimes(Q,\rho)\simeq(Q_1,\ba)\otimes(H_1,\ba).
\]
Therefore,
\begin{align*}
\Ad_{\pform{\lambda,\mu}}\otimes(Q,\rho) &\simeq
\Ad_{\qform{1,-\lambda}}\otimes (Q_1,\ba)\otimes(H_1,\ba)\\ &\in
\bigl((Q_1,\ba)\otimes(H_1,\ba)\bigr) \boxplus \bigl((Q_1,\ba)\otimes
(H_1,\ba)\bigr),
\end{align*}
and finally
\[
(A,\sigma)\in \bigl(Q_1,\ba)\otimes (H_1,\ba)\bigr)\boxplus
\bigl((Q_1,\ba)\otimes (H_1,\ba)\bigr) \boxplus
\bigl((M_2(F),\ba)\otimes (Q,\ba)\bigr).
\]
It follows that $\{0,[Q],[Q_1],[H_1]\}$ is a decomposition group for $(A,\sigma)$. 

If in addition $\sigma$ is hyperbolic, we may choose $\mu=1$, so that $\{[Q_1],[H_1]\}=\{0,[Q]\}$. 
Hence $\{0,[Q]\}$ is a decomposition group of $(A,\sigma)$ in this case. 

\end{ex} 

\begin{ex}
\label{ex:ind4}
If $A$ has index $4$, then all decomposition groups of $(A,\sigma)$ have order $8$. 
Indeed, since $[A]=[Q_i]+[H_i]$, the quaternion algebras
$Q_i$ and $H_i$ all are division algebras. Therefore $[Q_i]\not =0$ for $i=1,2,3$.  
Moreover, we have $[Q_1]+[Q_2]+[Q_3]=[A]$. Since $A$ has index $4$,
this guarantees $[Q_i]+[Q_j]$ is non zero if $i\not =j$. Therefore, $[Q_1],[Q_2],[Q_3]$ are $\mz/2$ linearly independent, and they do generate a group of order $8$. 

One may also check that the involution $\sigma$ is anisotropic in this case. 
Indeed, $A\simeq M_3(D)$ for some degree $4$ division algebra $D$,
hence $A$ does not carry any hyperbolic involution.  
Moreover, its isotropic involutions with trivial discriminant are
Witt-equivalent to $(Q_1,\ba)\otimes (H_1,\ba)$, for some quaternion
division algebras $Q_1$ and $H_1$ such that $D\simeq Q_1\otimes
H_1$. Hence isotropic involutions on $A$ with trivial discriminant
have non trivial Clifford invariant 
\[
[Q_1]+\{0,[D]\}=[H_1]+\{0,[D]\}\neq0\in H^2(F)/\{0,[D]\}.
\]

\end{ex}

\relax From these examples, we easily get the following
characterization of isotropy and hyperbolicity:  
\begin{lem}
\label{isodecgp.lem}
Let $(A,\sigma)$ be a degree $12$ algebra with orthogonal involution with trivial discriminant and trivial Clifford invariant. 
\begin{enumerate} 
\item[(i)]
The involution $\sigma$ is isotropic if and only if it admits a decomposition group generated by $[A]$ and $[Q_1]$ for some quaternion algebra $Q_1$.  
\item[(ii)] The involution $\sigma$ is hyperbolic if and only if it
  admits $\{0,[A]\}$ as a  decomposition group.  
\item[(iii)] The algebra with involution $(A,\sigma)$ is split and
  hyperbolic if and only if it admits $\{0\}$ as a decomposition
  group.  
\end{enumerate}

\end{lem} 

\begin{proof}
Assertion (iii) is clear from the definition of a decomposition group, since $(M_2(F),\ba)\otimes (M_2(F),\ba)$ is hyperbolic. 
For (i) and (ii), the direct implications immediately follow from the
previous examples.  
To prove the converse, let us first assume $(A,\sigma)$ admits $\{0,[A]\}$ as a decomposition group. 
Since this group has order $1$ or $2$, $A$ cannot have index $4$ by
Example~\ref{ex:ind4}. Therefore, it is Brauer-equivalent to a
quaternion algebra $Q$.  
Moreover, by definition, 
$$(A,\sigma)\in\boxplus_{i=1}^3 \bigl((Q,\ba)\otimes
(M_2(F),\ba)\bigr).$$ 
Since each summand is hyperbolic, this proves $\sigma$ is hyperbolic. 

Assume now that $U$ is generated by $[A]$ and $[Q_1]$. The
  order of $U$ then divides~$4$, hence by Example~\ref{ex:ind4} $A$ is
  Brauer-equivalent to some quaternion algebra $Q$. Thus,
$U=\{0,[Q],[Q_1],[H_1]\}$, with $M_2(Q)\simeq Q_1\otimes H_1$. If
$[Q_1]=0$, then the previous case applies, and $\sigma$ is
hyperbolic. Assume now $[Q_1]\not =0$.  
We get 
$$(A,\sigma)\in\boxplus_{i=1}^3(Q_i,\ba)\otimes (H_i,\ba),$$ with for
$i=2$, $3$ $\{[Q_i],[H_i]\}$ equal either to  
$\{0,[Q]\}$ or to $\{[Q_1],[H_1]\}$. Picking an arbitrary element in $\{[H_1],[Q_1]\}$ for $1\leq i\leq 3$, we get three quaternion algebras whose sum is never $0$. Therefore, since by Corollary~\ref{add-dec.cor} we have $[H_1]+[H_2]+[H_3]=0$, at least one summand must be $(Q,\ba)\otimes (M_2(F),\ba)$, and this proves $\sigma$ is isotropic. 
\end{proof} 

\begin{rem}
  \label{rem:anyU}
  Reversing the viewpoint, note that any subgroup $U\subset{}_2\br(F)$ of order~$8$ in which all the
  nonzero elements except at most one have index~$2$ is the
  decomposition group of some central simple algebra of degree~$12$
  with orthogonal involution of trivial discriminant and trivial
  Clifford invariant. If all the nonzero elements in $U$ have
  index~$2$, pick a quaternion algebra $D$ representing a nonzero
  element in $U$; otherwise, let $D$ be the division algebra such that
  $[D]\in U$ and $\ind D>2$. In each case, we may organize the other
  nonzero elements in $U$ in pairs $[Q_i]$, $[H_i]$ such that
  $[D]=[Q_i]+[H_i]$ for $i=1$, $2$, $3$, and
  $\sum_{i=1}^3[H_i]=0$. Any algebra with involution $(A,\sigma)$ in
  $\orthsum_{i=1}^3\bigl((Q_i,\ba)\otimes(H_i,\ba)\bigr)$ has
  decomposition group $U$. Modifying the scalars
  in the direct sum leads to several nonisomorphic such $(A,\sigma)$.
  Moreover, when all the nonzero elements in $U$ have index~$2$, we
  may select for $D$ various quaternion algebras, and thus obtain
  various $(A,\sigma)$ that are not Brauer-equivalent. Similarly, any
  subgroup 
  $U\subset{}_2\br(F)$ of order~$4$ containing at most one element
  $[D]$ with $\ind D>2$, and any subgroup $\{0,[Q]\}$ where $Q$ is a
  quaternion algebra, is the decomposition group of some central
  simple algebra of degree~$12$ endowed with an isotropic orthogonal
  involution of trivial 
  discriminant and trivial Clifford invariant.
\end{rem}

The decomposition groups of $(A,\sigma)$ are subgroups of the Brauer group generated by at most three quaternion algebras. 
Those subgroups were considered by Peyre in~\cite{Peyre}. His results will prove useful to study degree $12$ algebras with involution. For the reader's convenience, we recall them in the next section. 

\subsection{A complex of
  Peyre} 
\label{sec:Peyre}

Let $F$ be an arbitrary field, and let $U\subset\br F$ be a finite 
subgroup of the Brauer group of 
$F$. We let $F^\times \cdot U$ denote the subgroup of $H^3(F)$ generated
by classes $\lambda\cdot \alpha$, with $\lambda\in F^\times$ and
$\alpha\in U$; any element in $F^\times \cdot U$ can be written as
$\sum_{i=1}^r\lambda_i\cdot \alpha_i$ for some $\lambda_i\in
F^\times$, where $\alpha_1$,\ldots, $\alpha_r$ is a generating set for
the group $U$. Let $F_U$ be the function field of the product of the
Severi--Brauer varieties associated to elements of $U$. Clearly, $F_U$
splits all the elements of $U$, hence the subgroup $F^\times \cdot U$ vanishes
after scalar extension to $F_U$.  Therefore, the following sequence is a
complex, which was first introduced and studied by Peyre
in~\cite[\S4]{Peyre}: 
\[
F^\times \cdot U\ra H^3(F)\ra H^3(F_U).
\]
We let $\ch_U$ denote the corresponding homology group, that is
\[
\ch_U=\frac{\ker(H^3(F)\ra H^3(F_U))}{F^\times\cdot U}.
\]
We now return to our standing hypothesis that the characteristic of
$F$ is different from~$2$.
Peyre considers in particular subgroups $U\subset \br F$ generated by
the Brauer classes of at most three quaternion algebras, and
proves:

\begin{thm}[Peyre {\cite[Thm 5.1]{Peyre}}]
  \label{Peyre}
  If $U$ is generated by the Brauer classes of two quaternion
  algebras, then $\ch_U=0$.
\end{thm}

In the next section, we need to consider only subgroups $U$
such that all the elements of $U$ are quaternion algebras; we call
them \emph{quaternionic subgroups} of the Brauer group. These
subgroups have also been 
investigated by Sivatski \cite{Sivatski}. We have:

\begin{thm}[{\cite[Prop. 6.1]{Peyre}}, {\cite[Cor. 11]{Sivatski}}]
  \label{Peyre-cns}
  If $U\subset\br F$ is generated by the Brauer classes of three
  quaternion algebras, then $\ch_U=0$ or $\mz/2\mz$.
  Assume in addition $U$ is quaternionic. Then the following
  conditions are equivalent:
  \begin{enumerate}
    \item[(a)] $\ch_U=0$; 
    \item[(b)] $U$ is split by an extension of $F$ of degree~$2m$ for
      some odd $m$;
    \item[(c)] $U$ is split by a quadratic extension of $F$.
  \end{enumerate}
\end{thm}

The result that $\ch_U=0$ or $\mz/2\mz$ and the equivalence
(a)$\iff$(b) are due to Peyre~\cite[Prop. 6.1]{Peyre}. The equivalence
(b)$\iff$(c) was proved by Sivatski~\cite[Cor.~11]{Sivatski}.

We say that an extension $K$ of $F$ \emph{splits} a subgroup
  $U\subset \br F$ if it splits all the elements in $U$. If $K$ splits
  a decomposition group of a central simple algebra with orthogonal
  involution $(A,\sigma)$, then $A_K$ is
  split because $[A]\in U$, and $\sigma_K$ is hyperbolic by
  Lemma~\ref{isodecgp.lem}(iii) 
  because $(A_K,\sigma_K)$ has a trivial decomposition
  group. Therefore, Theorem~\ref{Peyre-cns} is relevant for the
  quadratic splitting of $(A,\sigma)$, as we will see in
  \S\,\ref{subsec:quadsplitanddec}.

\section{The Arason invariant and the homology of Peyre's complex}
\label{sec:ArasonPeyre}

As in the previous section, $(A,\sigma)$ is a degree $12$ algebra with orthogonal involution of trivial discriminant and trivial Clifford invariant. 
\relax From now on, we assume in addition that the Arason invariant
$e_3(\sigma)$ is well defined.  
So the algebra $A$ has even co-index, hence index $1$ or $2$. 
Under this assumption, any decomposition group of $(A,\sigma)$ is quaternionic, that is consists only of Brauer classes of quaternion algebras. 
In this section, we relate the decomposition groups of $(A,\sigma)$ with the values of the Arason invariant $e_3(\sigma)$. Reversing the viewpoint we then explain how one can use the Arason invariant to find explicit generators of the homology group $\ch_U$ of Peyre's complex, for any quaternionic subgroup $U\subset \br(F)$ of order dividing $8$. 

\subsection{Arason invariant in degree $12$.}

For orthogonal involutions on a degree $12$ algebra, isotropy and
hyperbolicity can be detected via the Arason invariant as follows:  
\begin{thm}
\label{iso.thm}
Let $(A,\sigma)$ be a degree $12$ and index $1$ or $2$ algebra with orthogonal involution of trivial discriminant and trivial Clifford invariant. 
\begin{enumerate}
\item[(i)] The involution $\sigma$ is hyperbolic if and only if $e_3(\sigma)=0\in M_A^3(F)$. 
\item[(ii)] The involution $\sigma$ is isotropic if and only if
    $e_3(\sigma)=e_3(\pi)+ F^\times\cdot[A]\in M_A^3(F)$ for some
    $3$-fold Pfister form $\pi$, i.e. $f_3(\sigma)=0$ and
    $e_3(\sigma)$ is represented by a symbol.
\end{enumerate} 
\end{thm}
 
\begin{proof} Assume first that $A$ is split, so that $\sigma$ is adjoint to a $12$-dimensional quadratic form $\varphi$, and $e_3(\sigma)=e_3(\varphi)\in {}_2H^3(F)$.   
Since the Arason invariant for quadratic forms has kernel the 4th
  power $I^4F$ of the fundamental ideal of the Witt ring $W(F)$, the
  first equivalence follows from the Arason--Pfister Hauptsatz.

  To prove (ii), note that there is no $10$-dimensional
 anisotropic quadratic form 
  in $I^3F$, see \cite[Prop.~XII.2.8]{Lam}.  So if $\varphi$
  is isotropic, then it has two hyperbolic 
  planes, and it is Witt-equivalent to a multiple of some $3$-fold
  Pfister form $\pi$. Hence, $e_3(\varphi)=e_3(\pi)$.  Assume
  conversely that $e_3(\varphi)=e_3(\pi)$.  By condition~(i),
  $\varphi$ becomes hyperbolic over the function field of
  $\pi$. Therefore, by~\cite[Th.~X.4.11]{Lam}, the
  anisotropic kernel of $\varphi$ is a multiple of $\pi$. In view of
  the dimensions, 
  this implies $\varphi=\qform{\alpha}\pi+2\mh$ for some $\alpha\in
  F^\times$. In particular, $\varphi$ is isotropic.
  
Assume now $A=M_6(Q)$ for some quaternion division algebra $Q$. 
  By a result of Dejaiffe \cite{Dej} and of
  Parimala--Sridharan--Suresh \cite[Prop 3.3]{PSS}, the involution
  $\sigma$ is hyperbolic if and only if it is hyperbolic after scalar
  extension to a generic splitting field $F_Q$ of the quaternion
  algebra $Q$. Since the restriction map $M_Q^3(F)\ra H^3(F_Q)$ is
  injective, the split case gives the result in index $2$.  
  If $\sigma$ is isotropic, its
  anisotropic part has degree $8$ and index $2$. The explicit description of the Arason invariant in
  degree~$8$ given 
  in~\cite[Th.~5.2]{QT:Arason} shows it is equal to $e_3(\pi)\mod
  F^\times\cdot [Q]\in M_Q^3(F)$ for some $3$-fold Pfister form $\pi$.   
  Assume conversely that
 $e_3(\sigma)=e_3(\pi)+ F^\times\cdot[Q]$. After scalar extension to $F_Q$, the split case shows $\sigma_{F_Q}$ is
  isotropic.  By Parimala--Sridharan--Suresh~\cite[Cor. 3.4]{PSS}, this
  implies $\sigma$ itself is isotropic.
\end{proof}

\begin{rem} 
\label{iso.rem}
In the split isotropic case, the involution can be explicitly
described from its Arason invariant: the proof of
  Theorem~\ref{iso.thm}(ii) shows that $\sigma$ is adjoint to $\pi+2\mh$
if $e_3(\sigma)=e_3(\pi)\in{}_2H^3(F)$.  
In index $2$, we also get an explicit description of
  $(A,\sigma)$ in the isotropic case. Indeed, 
we have
\[
(A,\sigma) = (M_2(Q),\hyp)\boxplus(M_4(Q),\sigma_0)
\]
for some orthogonal involution $\sigma_0$ with trivial discriminant
and trivial Clifford invariant, and $e_3(\sigma)=e_3(\sigma_0)$. If
$(a,b,c)$ is a symbol representing $e_3(\sigma)$, then
by~\cite[Th.~5.2]{QT:Arason} we may assume that one of the slots, say
$a$, is such that $F(\sqrt{a})$ splits $Q$, hence $Q$ carries an
orthogonal involution $\rho$ with discriminant~$a$. Theorem~5.2 of
\cite{QT:Arason} further shows that
\[
(M_4(Q),\sigma_0) \simeq (Q,\rho)\otimes \Ad_{\pform{b,c}}.
\]
\end{rem}

Under some additional condition, we also have the following classification result: 
\begin{prop}
  \label{qf.prop}
 Let $A=M_6(Q)$ be a degree $12$ algebra of index at most $2$, and let $\sigma$ and $\sigma'$ be two orthogonal involutions with trivial discriminant and trivial Clifford invariant. 
We assume either $A$ is split, or $\sigma$ is isotropic. 
The involutions $\sigma$ and $\sigma'$ are isomorphic if and only if $e_3(\sigma)=e_3(\sigma')$. 
  \end{prop}

\begin{proof}
It is already known that two isomorphic involutions have the same
Arason invariant, so we only need to prove the converse. Assume first
that $A$ has index $2$, in which case we assume in addition that
$\sigma$ is isotropic. Since $e_3(\sigma)=e_3(\sigma')$, by
Theorem~\ref{iso.thm}, the involution $\sigma'$ also is
isotropic. The result then follows from the explicit description given
in Remark~\ref{iso.rem}, or equivalently
from~\cite[Cor.~5.3(2)]{QT:Arason}, which shows that the anisotropic
parts of $\sigma$ and $\sigma'$ are isomorphic. 

   Assume now $A$ is split, and $\sigma$ and $\sigma'$ are adjoint to
   $\varphi$ and $\varphi'$ respectively. We have
   $e_3(\varphi)=e_3(\varphi')$. If there exists a $3$-fold 
  Pfister form $\pi$ such that $e_3(\varphi)=e_3(\varphi')=e_3(\pi)$,
  then $\varphi$ and $\varphi'$ are both similar to
  $\pi+2\mh$. Otherwise, they are anisotropic, and the result in this
  case follows by combining Pfister's theorem (see for
  instance~\cite[Th.~8.1.1]{Kahn}), which asserts that $\varphi$ and
  $\varphi'$ can be decomposed as tensor products of a
  $1$-fold Pfister form and an Albert form, with Hoffmann's
  result~\cite[Corollary]{Hoff}, which precisely says that two such
  forms are similar if and only if their difference is in $I^4F$.
  
  \end{proof}
  
\subsection{Arason invariant and decomposition groups}

Recall from Example~\ref{ex:split} that the decomposition groups
corresponding to additive decompositions of $(A,\sigma)$ are
quaternionic subgroups of order at most~$4$ when $A$ is split. Hence,
by Peyre's Theorem~\ref{Peyre}, the corresponding homology group is
trivial, $\ch_U=0$. Using this, we have:  

\begin{prop}
  Let $\varphi$ be a $12$-dimensional quadratic form in $I^3F$, and
  let $U=\{0,[Q_1],[Q_2],[Q_3]\}\subset\Br(F)$ be a quaternionic
  subgroup of order at most~$4$. For $i=1$, $2$, $3$, let $n_{i}$ be
  the norm form of $Q_i$. The following are
  equivalent:
  \begin{enumerate}
  \item[(a)] There exists $\alpha_1$, $\alpha_2$, $\alpha_3\in F^\times$
    such that $\varphi=\qform{\alpha_1}n_{1}\perp
    \qform{\alpha_2}n_{2}\perp \qform{\alpha_3}n_{3}$;
    \item[(b)] $U$ is a decomposition group of $\Ad_\varphi$;
  \item[(c)] $\varphi$ is hyperbolic over $F_U$;
  \item[(d)] $e_3(\varphi)\in \ker\bigl(H^3(F)\ra H^3(F_U)\bigr)$;
  \item[(e)] $e_3(\varphi)\in F^\times\cdot U$.
  \end{enumerate}
\end{prop}

\begin{proof}
  The equivalence between (a) and (b) follows from Example~\ref{split.ex}.
  Assume $\varphi$ decomposes as in (a). Since the field $F_U$ splits
  all three quaternion algebras $Q_i$, hence also their norm forms
  $n_{i}$, the form $\varphi$ is hyperbolic over $F_U$, hence assertions
  (c) and (d) hold. By Peyre's result~\ref{Peyre}, we also get (e),
  and it only remains to prove that (e) implies (a).

  Thus, assume now that $e_3(\varphi)\in F^\times\cdot U$. Since the
  subgroup $U$ is generated by $[Q_1]$ and $[Q_2]$, there exists
  $\lambda_1$ and $\lambda_2\in F^\times$ such
  that 
  \[
  e_3(\varphi)=(\lambda_1)\cdot [Q_1]+(\lambda_2)\cdot [Q_2].
  \]
  The product $Q_1\otimes Q_2\otimes Q_3$ is split, so by the common slot
  lemma (\cite[Th.~III.4.13]{Lam}), we may assume $Q_i=(a,b_i)_F$ for
  some $a$ and 
  $b_i\in F^\times$. A direct computation then shows that
  $n_1-n_2=\qform{b_2}n_3$. Hence the $12$-dimensional quadratic form
  $\qform{-\lambda_1}n_1+\qform{\lambda_2}n_2+\qform{b_2}n_3$ is Witt-equivalent to 
 $\qform{1,-\lambda_1}n_1+\qform{-1}\qform{1,-\lambda_2}n_2$, which 
  has the same Arason invariant as $\varphi$. By
  Proposition~\ref{qf.prop}, 
  this form
  is similar to $\varphi$, so that $\varphi$ has an additive
  decomposition as required.
\end{proof}

Let us consider now the index $2$ case. By Lemma~\ref{isodecgp.lem},
$(A,\sigma)$ admits decomposition groups of order $4$ if and only if
it is isotropic. We prove:  
\begin{prop}
\label{isotAra.prop}
Let $A=M_6(Q)$ be an algebra of index $\leq 2$, and consider an orthogonal involution $\sigma$ on $A$, with trivial discriminant and trivial Clifford invariant. Pick a subgroup $U=\{0,[Q],[Q_1],[H_1]\}\subset \br(F)$ containing the class of $Q$. The following are equivalent: 
\begin{enumerate}
\item [(a)] $(A,\sigma)$ admits an additive decomposition of the following type:
\[(A,\sigma)\in \bigl((Q_1,\ba)\otimes (H_1,\ba)\bigr)\boxplus\bigl((Q_1,\ba)\otimes(H_1,\ba)\bigr)\boxplus \bigl((M_2(F),\ba)\otimes (Q,\ba)\bigr);\]
\item[(b)] $U$ is a decomposition group of $(A,\sigma)$; 
\item[(c)] $\sigma$ is hyperbolic over $F_U$; 
\item[(d)] $e_3(\sigma)\in \ker\bigl(M_A^3(F)\ra H^3(F_U)\bigr)$; 
\item[(e)] There exists $\alpha\in F^\times$ such that $e_3(\sigma)=(\alpha)\cdot [Q_1]\mod F^\times\cdot [Q]\in M_Q^3(F)$. 
\end{enumerate} 
\end{prop} 

\begin{proof}
The proof follows the same line as for the previous proposition. 
By the definition of decomposition groups, (a) implies
(b). Conversely, if~(b) holds, then $(A,\sigma)$ has an
  additive decomposition with summands isomorphic to
  $(Q_1,\ba)\otimes(H_1,\ba)$ or to $(M_2(F),\ba)\otimes(Q,\ba)$. If
  $Q_1$ or $H_1$ is split, then the two kinds of summands are
  isomorphic, hence (a) holds. If $Q_1$ and $H_1$ are not split, then
  the number of summands isomorphic to $(Q_1,\ba)\otimes(H_1,\ba)$
  must be even because $e_2(\sigma)=0$ (see
  Proposition~\ref{prop:additinv}), and it must be nonzero because $U$
  is the corresponding decomposition group. Therefore, (a) holds.

Now, assume $(A,\sigma)$ satisfies (a). 
Since the field $F_U$ splits $Q$, $Q_1$ and $H_1$, (c)
holds. Assertion (d) follows since hyperbolic involutions have trivial
Arason invariant. 
By Peyre's Proposition~\ref{Peyre}, we deduce assertion (e), and it
only remains to prove that (e) implies (a).  
Hence, assume \[e_3(\sigma)=(\alpha)\cdot [Q_1]\mod F^\times\cdot [Q]\in M_Q^3(F),\] for some $\alpha\in F^\times$ and some quaternion algebra $Q_1$. 
By Theorem~\ref{iso.thm}(ii), the involution $\sigma$ is isotropic. 
Hence, in view of Proposition~\ref{qf.prop}, it is enough to find an involution $\sigma'$ satisfying $(a)$ and having $e_3(\sigma')=(\alpha)\cdot [Q_1]\mod F^\times [Q]$. 
Since $Q\otimes Q_1=H_1$ has index $2$, the quaternion algebras $Q$
and $Q_1$ have a common slot (see \cite[Th.~III.4.13]{Lam}).
Therefore, there exists $a$, $b$, $b_1\in F^\times$ such that $Q=(a,b)$ and
$Q_1=(a,b_1)$. Let $\rho$ be an orthogonal involution on $Q$ with
discriminant $a$, and let
\[(A,\sigma')=(M_2(Q),{\mathrm {hyp}})\orthsum \bigl((Q,\rho)\otimes
  \Ad_\pform{b_1,\alpha}\bigr).
\] 
One component of the Clifford algebra of 
$(Q,\rho)\otimes \Ad_\pform{b_1}$ is given by the cup product of the
discriminants of $\rho$ and $\ad_\pform{b_1}$, that is
$(a,b_1)=Q_1$. By Remark~\ref{iso.rem}, it follows that
$e_3(\sigma')=e_3(\sigma)$, hence 
$(A,\sigma)\simeq(A,\sigma')$ by Proposition~\ref{qf.prop}. Therefore
$(A,\sigma)$ satisfies (a) as required. 
\end{proof}

\subsection{Generators of the homology $\ch_U$ of Peyre's complex}

Let $A=M_6(Q)$ for some quaternion $F$-algebra $Q$, and let
  $\sigma$ be an orthogonal involution on $A$ with trivial
  discriminant and trivial Clifford invariant. Consider an additive
decomposition of $(A,\sigma)$ as in Theorem~\ref{thm:adddec12},
\[
(A,\sigma) \in \orthsum_{i=1}^3\bigl((Q_i,\ba)\otimes (H_i,\ba)\bigr),
\] 
and let $U$ be the corresponding decomposition group, which is a
quaternionic subgroup of $\Br(F)$,
\[
U=\{0,[Q],[Q_1],[H_1],[Q_2],[H_2],[Q_3],[H_3]\}.
\]
Since $F^\times\cdot [Q]\subset F^\times\cdot U$, we may consider the
canonical map
\[
\ba^U\colon M^3_Q(F)\to H^3(F)/F^\times\cdot U.
\]
As in \S\ref{sec:Peyre}, let $F_U$ be the function field of the product of the
Severi--Brauer varieties associated to elements of $U$. Since $F_U$
splits $U$, Lemma~\ref{isodecgp.lem} shows that $A_{F_U}$ is split
and $\sigma_{F_U}$ is hyperbolic, hence
$e_3(\sigma)_{F_U}=0$. Therefore, $\overline{e_3(\sigma)}^U$ lies in
the homology $\ch_U$ of Peyre's complex. 

As explained in Remark~\ref{rem:anyU}, for any quaternionic subgroup
$U\subset \br(F)$ of order dividing $8$, we may find algebras with
involution $(A,\sigma)$ for which $U$ is a decomposition group.  
The main result of this section is:

\begin{thm}
  \label{hom.thm}
  Let $U$ be a quaternionic subgroup of $\br F$ of order dividing $8$.
  For any $(A,\sigma)$ admitting $U$ as a decomposition group, the
  class of the Arason invariant $\overline{e_3(\sigma)}^U$ is a
  generator of the homology group $\ch_U$ of Peyre's
  complex. 
  \end{thm}

The main tool in the proof is the following proposition: 
\begin{prop}
  \label{same-dec.prop}
  Let $U$ be a quaternionic subgroup of $\br F$ of order dividing $8$, and pick an algebra $A=M_6(Q)$ with orthogonal involution $\sigma$, admitting 
  $U$ as a decomposition group.
  \begin{enumerate}
  \item[(i)]
  For all involutions $\sigma'$ on $A$ such that $(A,\sigma')$ also admits $U$ as a decomposition group, we have 
  \[
  \overline{e_3(\sigma')}^U=\overline{e_3(\sigma)}^U.
  \]
  \item[(ii)]
  Conversely, for all $\xi\in M^3_A(F)$ such that $\overline{\xi}^U =
  \overline{e_3(\sigma)}^U$, there exists an involution
  $\sigma'$ on $A$ such that $U$ is a decomposition group of
  $(A,\sigma')$ and $e_3(\sigma')=\xi\mod F^\times\cdot[A]$.
  
 \item[(iii)] There exists a hyperbolic involution $\sigma'$ on $A$ admitting $U$ as a decomposition group 
 if and only if $\overline{e_3(\sigma)}^U=0$. 
 
  \end{enumerate}
\end{prop}

\begin{proof}
  (i)
  Since $U$ is a decomposition group of $(A,\sigma)$ and
  $(A,\sigma')$, we have
  \[
  (A,\sigma)\mbox{ and }(A,\sigma')\in
  \orthsum_{i=1}^3\bigl((Q_i,\ba)\otimes (H_i,\ba)\bigr).
  \] 
  Therefore, $\sigma$
  and $\sigma'$ are adjoint to some skew-hermitian forms $h$ and $h'$
  over $(Q,\ba)$
  satisfying \[h=h_1\perp h_2\perp h_3\mbox{ and }h'=\qform{\alpha_1}h_1\perp\qform{\alpha_2}h_2\perp\qform{\alpha_3}h_3,\]
  for some $h_i$ such that $\ad_{h_i}\simeq
  (Q_i,\ba)\otimes (H_i,\ba)$, and some $\alpha_i\in F^\times$.  Therefore,
  \[
  e_3(\sigma)-e_3(\sigma')=e_3(\perp_{i=1}^3\qform{1,-\alpha_i}h_i).
  \]
  Since $h_i$ has discriminant $1$, Proposition~\ref{compute.prop}
  applies to each summand and shows $\qform{1,-\alpha_i}\otimes h_i$
  has trivial discriminant and trivial Clifford invariant, and
  \[
  e_3(\qform{1,-\alpha_i}\otimes h_i)=\alpha_i\cdot
  e_2(h_i)=\alpha_i\cdot [Q_i]\in M_Q^3(F).
  \]
  Therefore, $e_3(\sigma)-e_3(\sigma')$ is represented modulo
  $F^\times\cdot\{0,[Q]\}$ by $\sum_{i=1}^3\alpha_i\cdot
  [Q_i]$. Since this element lies in $F^\times\cdot U$, we have
  $\overline{e_3(\sigma)}^U = \overline{e_3(\sigma')}^U$.

  (ii)
  Consider a skew hermitian form $h$ over $(Q,\ba)$ such that $\sigma=\ad_h$, and a decomposition $h=h_1\perp h_2\perp h_3$  as in the proof of (i).  Since $\overline{\xi}^U =
  \overline{e_3(\sigma)}^U$, the difference $e_3(\sigma)-\xi\in
  M^3_Q(F)$ is represented by a cohomology class of the form
  $\sum_{i=1}^3\alpha_i\cdot[Q_i]$ for some $\alpha_i\in
  F^\times$. The computation in~(i) shows that $e_3(\ad_{h'})=\xi$ for
  $h'=\qform{\alpha_1}h_1\perp\qform{\alpha_2}h_2\perp\qform{\alpha_3}h_3$. 

(iii) 
  It follows from (ii) that $\overline{e_3(\sigma)}^U=0$
  if and only if there exists an involution $\sigma'$ with
  decomposition group $U$ and $e_3(\sigma')=0$.  
  Theorem~\ref{iso.thm}(i) completes
  the proof by showing $\sigma'$ is hyperbolic.
\end{proof}

With this in hand, we can now prove Theorem~\ref{hom.thm}.

\begin{proof}[Proof of Theorem~\ref{hom.thm}]
  Since $\ch_U$ is either $0$ or $\mz/2\mz$, 
 in order to prove that
  $\overline{e_3(\sigma)}^U$ generates $\ch_U$ it is enough to prove
  that $\ch_U$ is trivial as soon as $\overline{e_3(\sigma)}^U=0$. If
  $U$ has order at most $4$, then $\ch_U$ is trivial
  by Theorem~\ref{Peyre}. Hence, let us assume $U$ has order $8$, and
  $\overline{e_3(\sigma)}^U=0$. By Proposition~\ref{same-dec.prop}, replacing $\sigma$ by $\sigma'$, we may assume 
  $\sigma$ is hyperbolic. Recall $\sigma$ is adjoint to a
  skew-hermitian form $h$, which admits a decomposition
  $h=h_1\perp h_2\perp h_3$ with $\Ad_{h_i}=(Q_i,\ba)\otimes (H_i,\ba)$. Since
  $U$ has order $8$, each summand $h_i$ is anisotropic. The
  hyperbolicity of $h$ says $h_1\perp h_2\simeq-h_3\perp\mh$ is
  isotropic. Therefore, there exists a pure quaternion $q$ such that
  $h_1$ and $h_2$ represent $q$ and $-q$ respectively. Over the
  quadratic extension $F(q)$ of $F$, the involutions $\ad_{h_1}$ and
  $\ad_{h_2}$ are isotropic. Since they are adjoint to $2$-fold
  Pfister forms, they are hyperbolic. Hence $F(q)$ splits the Clifford
  algebra of $h_1$ and $h_2$, that is the quaternion algebras $Q_1$,
  $Q_2$. Since the Brauer classes of $Q$, $Q_1$ and $Q_2$
  generate $U$, it follows that $F(q)$ is a quadratic splitting
  field of $U$. By Peyre's Theorem~\ref{Peyre-cns}, we get $\ch_U=0$ as
  required.
\end{proof}

\section{Quadratic splitting and the $f_3$ invariant}
\label{sec:quadsplit}

The $f_3$ invariant of an involution $\sigma$ vanishes if the underlying algebra $A$ is split, or of degree $\leq 10$. 
We keep focusing on the case of
  degree~$12$ algebras, where we have explicit examples with
$f_3(\sigma)\not=0$, see Corollary~\ref{f3.cor}. 
Thus, as in \S\,\ref{sec:ArasonPeyre}, $(A,\sigma)$ is a degree $12$
algebra with 
orthogonal involution of trivial discriminant and trivial Clifford
invariant for which the Arason and the $f_3$ invariants are
defined. In particular, $A$ has index at most $2$. 

Our first goal
is to characterize the vanishing of $f_3(\sigma)$; this is done in
Proposition~\ref{f3zero.prop} below.
As pointed out in Proposition~\ref{quadsplit.prop}, $f_3(\sigma)$ vanishes
if there exists a quadratic extension $K/F$ over which $(A,\sigma)$ is
split and hyperbolic.  Note that since $A$ is
Brauer-equivalent to a quaternion algebra,
there exist quadratic extensions of the base field $F$ over which $A$
is split. Moreover, using the additive decompositions of
Corollary~\ref{add-dec.cor}, one may easily find quadratic extensions
of the base field over which the involution is hyperbolic: it
suffices to consider a common subfield of the quaternion algebras
$H_1$, $H_2$, $H_3$, which exists by~\cite[Th.~III.4.13]{Lam} since
$[H_1]+[H_2]+[H_3]=0$. Yet, we give in Corollary~\ref{cor:comeback}
examples showing that the
converse of Proposition~\ref{quadsplit.prop} does not hold in degree
$12$: we may have $f_3(\sigma)=0$ even when there is no quadratic
extension that simultaneously splits $A$ and makes $\sigma$ hyperbolic.

First, we use quadratic forms to introduce an invariant of
quaternionic subgroups of the Brauer group of $F$, which, as we next
prove, coincides with the $f_3$-invariant of involutions admitting
this subgroup as a decomposition group.  

\subsection{The invariants $f_3(U)$ and $f_3(\sigma)$}

To any quaternionic subgroup $U$ of $\br(F)$, we may associate in a natural way a
quadratic form $n_U$ by taking the sum of the norm forms
$n_H$ of the quaternion algebras $H$ with Brauer class in $U$. We have:

\begin{lem} 
  \label{lem:f3U}
  Let $U$ be a quaternionic subgroup of $\br F$ generated by the
  Brauer classes of three quaternion algebras. The quadratic form
    $n_U=\sum_{[H]\in U} n_H$ satisfies $n_U\in I^3F$.
\end{lem}

\begin{proof}
  Pick three generators $[Q_1]$, $[Q_2]$ and $[Q_3]$ of $U$, and let
  $H_1$, $H_2$, $H_3$, $Q$ be quaternion algebras with Brauer classes
  $[H_1]=[Q_2]+[Q_3]$, $[H_2]=[Q_1]+[Q_3]$, $[H_3]=[Q_1]+[Q_2]$, and
  $[Q]=[Q_1]+[Q_2]+[Q_3]$. We have $[H_1]+[H_2]+[H_3]=0$, and 
  \[
  U=\{0, [Q],[Q_1],[H_1],[Q_2],[H_2],[Q_3],[H_3]\}.
  \]
  Since the difference $n_{Q_i}-n_{H_i}$ is Witt-equivalent to an
  Albert form of $Q_i\otimes H_i$, which is Brauer-equivalent to $Q$,
  there exists $\lambda_i\in F^\times$ such that in the Witt group of
  $F$, we have 
  $n_{Q_i}-n_{H_i}=\qform{\lambda_i}n_Q\in WF$. Therefore,
  \begin{equation}
  \label{eq:nU}
  n_U=\qform{1,\lambda_1,\lambda_2,\lambda_3}n_Q+
  \qform{1,1}(n_{H_1}+n_{H_2}+n_{H_3}).
  \end{equation}
Since the right side is in $I^3F$, the lemma is proved.
\end{proof}

In view of Lemma~\ref{lem:f3U}, we may associate to $U$ a cohomology
class of degree~$3$ as follows:

\begin{defi}
  For any quaternionic subgroup $U$ generated by three elements, we
  let $f_3(U)$ be the Arason invariant of the quadratic form $n_U$:
  \[
  f_3(U)=e_3(n_U)\in {}_2H^3(F).
  \]
\end{defi}

We may easily compute $f_3(U)$ from formula~\eqref{eq:nU}: 
Since $[H_1]+[H_2]+[H_3]=0$, we have $n_{H_1}+n_{H_2}+n_{H_3}\in
I^3F$, hence $\qform{1,1}(n_{H_1}+n_{H_2}+n_{H_3})\in I^4F$ and therefore
\begin{equation}
  \label{f3.eq}
  f_3(U)=(\lambda_1\lambda_2\lambda_3)\cdot [Q].
\end{equation}

With this in hand, we get:

\begin{prop}
  \label{hufu.prop}
  If $\ch_U=0$, then $f_3(U)=0$.
\end{prop}

\begin{proof}
  By Theorem~\ref{Peyre-cns}, if $\ch_U=0$ then $U$ admits
  a quadratic splitting field, i.e. the generators of $U$ have a
  common quadratic subfield. So there exist $a$, $b_1$, $b_2$, and
  $b_3\in 
  F^\times$ such that $Q_i=(a,b_i)_F$ for $i=1$, $2$, $3$. Thus, we have
  $H_1=(a,b_2b_3)_F$ and
  \[
  n_{Q_1}-n_{H_1}=\pform{a}(\pform{b_1}-\pform{b_2b_3})=
  \pform{a}\qform{-b_1,b_2b_3}=\qform{-b_1}n_Q.
  \] 
  Similar formulas hold for $i=2$, $3$, and we get
  \[
  f_3(U)=(-b_1b_2b_3)\cdot Q=(-b_1b_2b_3,a,b_1b_2b_3)=0\in
  {}_2H^3(F).
  \]
\end{proof}

In~\cite{Sivatski}, Sivatski asks about the
converse\footnote{Sivatski's invariant has a different definition, but
  one may easily check the quadratic form he considers is equivalent
  to $n_U$ modulo $I^4 F$.}, that is: if $f_3(U)=0$, does the homology
group $\ch_U$ vanish, or equivalently by Peyre's Theorem~\ref{Peyre-cns},
do the generators $Q_1$, $Q_2$, and $Q_3$ of the group $U$ have a
common quadratic subfield? Corollary~\ref{decsub.cor} below shows that
this is not the case. 

The relation between $f_3(U)$ and the $f_3$-invariant for involutions is given by the following: 
\begin{thm}
\label{f3.thm}
Let $(A,\sigma)$ be a central simple algebra of degree $12$ and index $\leq 2$, with orthogonal involution of trivial discriminant and trivial Clifford invariant. Let $U$ be a quaternionic subgroup of the Brauer group, generated by three elements. If $U$ is a decomposition group for $(A,\sigma)$ then $f_3(\sigma)=f_3(U)$. 
\end{thm} 

\begin{rem}
\label{order.rem}
(i) It follows that any two decomposition groups of a given algebra with involution have the same $f_3$-invariant, and any two algebras with involution having $U$ as a decomposition group have the same $f_3$-invariant. 

(ii) Let $c$ be a generator of $\ch_U$, and pick an arbitrary $(A,\sigma)$ having $U$ as a decomposition group. 
In view of theorem~\ref{hom.thm}, we have $e_3(\sigma)=c\mod F^\times\cdot U$. Hence $f_3(U)=f_3(\sigma)=2c\in {}_2H^3(F)$. 
In particular, $f_3(U)=0$ if and only if the homology group $\ch_U$ is generated by cohomology class of order $2$. 
\end{rem}

\begin{proof}[Proof of Theorem~\ref{f3.thm}]
The result follows from the computation of
  $f_3(\sigma)$ in Proposition~\ref{prop:f3sum} and the computation of
  $f_3(U)$ in \eqref{f3.eq}. 
  We use the same notation as in Definition~\ref{defi:decsubgrp}, and we let $h_i$ be a rank $2$ skew-hermitian form over $(Q,\ba)$ such that \[\Ad_{h_i}\simeq(Q_i,\ba)\otimes (H_i,\ba)\mbox{ and 
  }\sigma=\ad_{h_1\perp h_2\perp h_3}.\] 
  For $i=1$, $2$, $3$, let $q_i\in Q$ be a
  nonzero pure quaternion represented by $h_i$, and let $a_i=q^2_i\in
  F^\times$. Let also $b_i\in F^\times$ be such that
  $Q=(a_i,b_i)_F$. Scalar extension to $F(q_i)$ makes $h_i$ isotropic,
  hence hyperbolic since the discriminant of $h_i$ is
  trivial. Therefore, we have $h_i\simeq\qform{q_i}\qform{1,-\lambda_i}$ for
  some $\lambda_i\in F^\times$. The two components of the Clifford
  algebra of $\Ad_{h_i}$ are $(a_i,\lambda_i)_F$ and
  $(a_i,\lambda_ib_i)_F$, therefore
  \[
  \{Q_i,H_i\} = \{(a_i,\lambda_i)_F, (a_i,\lambda_ib_i)_F\}
  \qquad\text{for $i=1$, $2$, $3$.}
  \]
Since $Q$ contains a pure quaternion which anticommutes with $q_i$ and with square $b_i$, the form $h_i$ is isomorphic to $\qform{q_i}\qform{1,-\lambda_ib_i}$ for $i=1$, $2$,
  $3$. Replacing some $\lambda_i$ by $\lambda_ib_i$ if necessary, we may assume  $H_i=(a_i,\lambda_i)_F$ for all~$i$. 
  Since $[H_1]+[H_2]+[H_3]=0$, we get 
  $\sum_{i=1}^3(a_i,\lambda_i)_F=0$. By Proposition~\ref{prop:f3sum} this implies $f_3(\sigma)=\lambda_1\lambda_2\lambda_3\cdot[Q]$. On the other
  hand, since $n_{Q_i}-n_{H_i}=\pform{a_i,\lambda_ib_i}-\pform{a_i,\lambda_i}=\qform{\lambda_i}n_Q$, we have $f_3(U)=\lambda_1\lambda_2\lambda_2\cdot[Q]$
  by~\eqref{f3.eq}.
\end{proof}

\subsection{Quadratic splitting, the $f_3$ invariant, and
  decomposition groups} 
\label{subsec:quadsplitanddec}

By using Theorem~\ref{f3.thm} and Peyre's Theorem~\ref{Peyre-cns},
we can now translate in terms of
decomposition groups the two conditions we want to compare, as
follows:

\begin{prop}
\label{f3zero.prop}
  Let $(A,\sigma)$ be a degree $12$ and index $\leq 2$ algebra with orthogonal involution of trivial discriminant and trivial Clifford invariant. 
  The following conditions are equivalent: 
  \begin{enumerate}
  \item[(a)]
    $f_3(\sigma)=0$;
  \item[(b)]
    $(A,\sigma)$ has a decomposition group $U$ with $f_3(U)=0$;
  \item[(c)]
    $f_3(U)=0$ for all decomposition groups $U$ of $(A,\sigma)$.
  \end{enumerate}
  Likewise, the following conditions are equivalent:
      \begin{enumerate}
  \item[(a')]
    there exists a quadratic extension $K$ of $F$ such that $A_K$ is
    split and $\sigma_K$ is hyperbolic;
  \item[(b')]
    $(A,\sigma)$ has a decomposition
    group $U$ with $\ch_U=0$.
  \end{enumerate}
  Moreover, any of the conditions~\emph{(a')}, \emph{(b')} implies the
  equivalent conditions~\emph{(a)}, \emph{(b)}, \emph{(c)}.
\end{prop}

\begin{proof} 
The equivalence between conditions ~{(a)}, {(b)}, {(c)} follows
directly from Theorem~\ref{f3.thm}. Moreover, they can be deduced
from~{(a')}, {(b')} by Proposition~\ref{hufu.prop} or Proposition~\ref{quadsplit.prop}. 
Hence, it only
remains to prove that~{(a')} and {(b')} are equivalent. 

Assume first that $(A,\sigma)$ has a decomposition group $U$ with $\ch_U=0$. 
By Peyre's characterization of the vanishing of $\ch_U$ for
quaternionic groups, recalled in Theorem~\ref{Peyre-cns}, $U$ is split
by a quadratic extension $K$ of $F$. Hence, $(A_K,\sigma_K)$ admits
$\{0\}$ as a decomposition group. By Lemma~\ref{isodecgp.lem}, this
implies $(A_K,\sigma_K)$ is split and hyperbolic.  

To prove the converse, let us assume there exists a quadratic field extension 
  $K=F(d)$, with $d^2=\delta\in F^\times$, such that $A_K$ is split and
  $\sigma_K$ is hyperbolic.
  If $A$ is split, as explained in example~\ref{ex:split}, all
  decomposition subgroups $U$ of $(A,\sigma)$ have order dividing $4$,
  and therefore satisfy $\ch_U=0$ by Peyre's Theorem~\ref{Peyre}. 
 Assume next $\ind A=2$. Since $A_K$ is split, we may identify
  $K=F(d)$ with a subfield of the quaternion division algebra $Q$
  Brauer-equivalent to $A$, and thus consider $d$ as a pure quaternion
  in $Q$. Let $h$ be a skew-hermitian form over $(Q,\ba)$ such that
  $\sigma=\ad_h$. Since $h_K$ is hyperbolic, it follows from
  \cite[Prop., p.~382]{QT} that
  $h\simeq \qform d \varphi_0$ for some $6$-dimensional quadratic form
  $\varphi_0$ over $F$. Decompose  
  \[
  \varphi_0=\qform{\alpha_1}\qform{1,-\beta_1} \perp
  \qform{\alpha_2}\qform{1,-\beta_2} \perp
  \qform{\alpha_3}\qform{1,-\beta_3}\qquad\text{for some $\alpha_i$, $\beta_i\in F^\times$},
  \]
  and let $Q_i=(\delta,\beta_i)_F$ be the quaternion $F$-algebra with
  norm $n_{Q_i}=\pform{\delta,\beta_i}$ for $i=1$, $2$, $3$.
  Computation shows that
  $e_2\bigl(\qform{\alpha_id}\qform{1,-\beta_i}\bigr)$ is represented
  by $Q_i$ in $M^2_Q(F)$, hence $(A,\sigma)$ decomposes as 
  \[
  (A,\sigma)\in \orthsum_{i=1}^3\Ad_{\qform{\alpha_id}\qform{1,-\beta_i}}. 
  \]
 So, the subgroup $U\subset\br F$ generated by $[Q_1]$, $[Q_2]$ and
  $[Q_3]$ is a decomposition group for $(A,\sigma)$. Again, $U$ is
  split by $K$, hence $\ch_U=0$.  
\end{proof}

\subsection{Trivial $f_3$-invariant without quadratic splitting}
\label{ex.sec}

We now construct an algebra with involution $(A,\sigma)$, of degree
$12$ and index $2$, such that $f_3(\sigma)=0$, and yet, there is no
quadratic extension $K$ of $F$ over which $(A,\sigma)$ is both split
and hyperbolic. In particular, by Peyre's Theorem~\ref{Peyre-cns}, we
have $\ch_U\neq 0$ for all decomposition groups $U$ of $(A,\sigma)$.  
(See Remark~\ref{rem:weaker} for an example where $(A,\sigma)$ has a
decomposition group $U$ whith $\ch_U\neq0$ and another $U'$ with
$\ch_{U'}=0$.)

\begin{rem}
In his paper \cite[\S6.2]{Peyre}, Peyre provides an example of a
quaternionic subgroup $U\subset\Br(F)$ with $\ch_U\neq0$, but the way
he proves $\ch_U$ is nonzero is by describing an element $c\in H^3(F)$
which is not of order $2$, hence does not belong to $F^\times\cdot U$,
and yet is in the kernel of the restriction map $H^3(F)\ra
H^3(F_U)$. Thus, the group $U$ in Peyre's example satisfies
$f_3(U)\neq0$ (see remark~\ref{order.rem}). In this section, we construct an example of a different
flavor, namely a subgroup $U$ with $\ch_U\neq0$, but
$f_3(U)=0$. Hence, the homology group in this case is generated by a
cohomology class which is of order $2$, and in the kernel of
$H^3(F)\ra H^3(F_U)$, but does not belong to $F^\times \cdot U$.
\end{rem}

\begin{notation}
  \label{not:not}
  Until the end of this section, $k$ is a field (of characteristic
  different from~$2$), $M$ is a triquadratic field extension of $k$
  (of degree~$8$) and $K$ is a quadratic extension of $k$ in $M$,
  \[
  M=k(\sqrt a,\sqrt b,\sqrt c)\supset K=k(\sqrt a).
  \]
  We let $C$ be a central simple $k$-algebra of exponent~$2$ split by
  $M$ and we write
  \[
  [C]\in\Dec(M/k)
  \]
  to express the property that there exist $\alpha$, $\beta$,
    $\gamma\in k^\times$ such that
  \[
  [C]=(a,\alpha)_k+(b,\beta)_k+(c,\gamma)_k.
  \]
  
  The existence of algebras $C$ as above such that
  $[C]\notin\Dec(M/k)$ is shown in~\cite[\S5]{ELTW}. By contrast, it
  follows from a 
  theorem of Albert that every central simple algebra of exponent~$2$
  split by a biquadratic extension has a decomposition up to
  Brauer-equivalence into a tensor product of quaternion algebras
  adapted to the biquadratic extension (see \cite[Prop.~5.2]{LLT}), so
  (viewing $M$ as 
  $K(\sqrt{bc},\sqrt c)$) there exist $x$,
  $y\in K^\times$ such that
  \[
  [C_K]=(bc,x)_K+(c,y)_K.
  \]
  By multiplying $x$ and $y$ by squares in $K$, we may---and
  will---assume $x$, $y\notin k$.
  We have $\cores_{K/k}[C_K]=2[C]=0$, hence letting $N$ denote the
  norm map from $K$ to $k$, we obtain from the previous equation by
  the projection formula: $(bc,N(x))_k+(c,N(y))_k=0$. We may then
  consider the following quaternion $k$-algebra:
  \begin{equation}
    \label{quat.eq}
    H=(bc,N(x))_k=(c,N(y))_k.
  \end{equation}
  Since $N(x)$, $N(y)$ are norms from $K=k(\sqrt a)$ to $k$, we have
  $(a,N(x))_k=(a,N(y))_k=0$, hence we may also write
  \begin{equation}
    \label{quat.eq2}
    H=(abc,N(x))_k=(ac,N(y))_k.
  \end{equation}
  Let $B=(bc,x)_K\otimes_K(c,y)_K$ be the biquaternion algebra Brauer-equivalent to $C_K$, and let $\psi$ be the Albert form
  of $B$ over $K$ defined by
  \[
  \psi=\qform{bc,x,-bcx,-c,-y,cy}.
  \]
  Let $s\colon\,K\ra k$ be a nontrivial linear map such that $s(1)=0$,
  and let $s_\star$ denote the corresponding Scharlau transfer. Using
  the properties of $s_\star$ (see for instance~\cite[p.~189,
  p.~198]{Lam}), 
  we can make the following computation in the Witt group $W(k)$:
  \begin{align*}
    s_\star(\psi)=& s_\star(\qform{x,-bcx,-y,cy})=
    s_\star\bigl(\qform{x}\bigr)\pform{bc}-s_\star\bigl(\qform{y}\bigr)\pform{c}\\
    =& \qform{s(x)}\pform{bc,N(x)}-\qform{s(y)}\pform{c,N(y)}.
  \end{align*}
  (Recall that we assume $x$, $y\notin k$, so $s(x)$, $s(y)\neq0$.)
  In view of~(\ref{quat.eq}), the last equation yields
  \[
  s_\star(\psi)=\qform{s(x),-s(y)}n_H,
  \]
  where $n_H$ is the norm form of $H$. Thus, $s_\star(\psi)\in I^3(k)$, and we may consider
  \begin{equation}
    \label{eq:e3psi}
    e_3(s_\star(\psi))=s(x)s(y)\cdot[H]\in {}_2H^3(k).
  \end{equation}
  This class represents an invariant
  of $B$ defined by Barry~\cite{Barry}. It is shown in
  \cite[Prop.~4.4]{Barry} that 
  $e_3(s_\star(\psi))\in N(K^\times)\cdot[C]$ if and only if the
  biquaternion algebra $B$ has a descent to $k$, i.e., there exist
  quaternion $k$-algebras $A_1$, $A_2$ such that
  $B\simeq A_1\otimes_kA_2\otimes_kK$.

  Finally, let $t$ be an indeterminate over $k$, and let
  $F=k(t)$. Consider the subgroup $U\subset \br(F)$ generated by 
  the Brauer classes $(a,t)_{F}$,
  $(b,t)_{F}$ and $(c,t)_{F}+[H_{F}]$. In view of \eqref{quat.eq} and
  \eqref{quat.eq2}, 
  one may easily check that $U$ is a quaternionic subgroup of order
  $8$:
  \begin{multline*}
  U=\{0,(a,t)_F,(b,t)_F,(c,N(y)t)_F,\\ (ab,t)_F, (ac,N(y)t)_F,
  (bc,N(x)t)_F, (abc,N(x)t)_F\}.
  \end{multline*}
  We set
  \[
  \xi=t\cdot [C] + e_3(s_\star(\psi))\in{}_2H^3(F).
  \]
\end{notation}

This construction yields examples with trivial $f_3$ but with no
quadratic splitting mentioned in the introduction to
this section, as we proceed to show. First, we describe the group
$\ch_U$ and give a criterion for its vanishing:
 
\begin{thm}
  \label{thm:xi}
  Use the notation~\ref{not:not}. Denote by $\overline{\xi}^U\in
  H^3(F)/F^\times\cdot U$ 
  the image of $\xi\in H^3(F)$. We have
  \[
  \ch_U=\{0,\overline{\xi}^U\}\qquad\text{and}\qquad f_3(U)=0.
  \]
  Moreover, the following conditions are equivalent:
  \begin{enumerate}
  \item[(a)]
    $[C]\in\Dec(M/k)$;
  \item[(b)]
  $U$ is split by some quadratic field extension $E/F$; 
  \item[(c)]
    $\ch_U=0$;
  \item[(d)]
    $\xi\in F^\times\cdot U$. 
  \end{enumerate}
\end{thm}

The core of the proof is the following technical lemma:
 
\begin{lem}
  \label{lem:tech}
  With the notation~\ref{not:not}, every field extension of $F$ that splits $U$ also splits $\xi$.
\end{lem}

\begin{proof}
  Let $L$ be an extension of $F$ that splits $U$. We consider two
  cases, depending on whether $a\in L^{\times2}$ or $a\notin
  L^{\times2}$. Suppose first $a\in L^{\times2}$, so we may identify
  $K$ with a subfield of $L$, hence $x$, $y\in L^\times$ and
  $[C_L]=(bc,x)_L+(c,y)_L$. Since $L$ splits $(b,t)_F$, we have
  $(t,bc,x)_L=(t,c,x)_L$, hence $t\cdot[C_L]=xy\cdot(t,c)_L$. Since
  $L$ also splits $(t,c)_F+[H_F]$, we have
  \[
  t\cdot[C_L]=xy\cdot[H_L].
  \]
  Comparing with~\eqref{eq:e3psi}, we see that it suffices to show
  $xy\cdot[H_L]=s(x)s(y)\cdot[H_L]$ to prove that $L$ splits
  $\xi$.

  Let $\iota$ be the nontrivial automorphism of $K$ over $k$. Writing
  $x=x_0+x_1\sqrt a$ and $y=y_0+y_1\sqrt a$ with $x_i,y_i\in k$, we have 
  \[s(x)s(y)=x_1y_1s(\sqrt a)^2\mbox{ and }
  (x-\iota(x))(y-\iota(y))=4x_1y_1a.\] Hence $s(x)s(y)\equiv
  (x-\iota(x))(y-\iota(y))\bmod L^{\times2}$. We also have
  \[(x-\iota(x),N(x))_K = (x,N(x))_K\mbox{ because }
  (x^2-N(x),N(x))_K=0.\] From the expression $H=(bc,N(x))_k$ it then
  follows that $x\cdot[H_K]=(x-\iota(x))\cdot[H_K]$. Similarly, from
  $H=(c,N(y))_k$ we have $y\cdot[H_K]=(y-\iota(y))\cdot[H_K]$, hence
  \[
  xy\cdot[H_L]=s(x)s(y)\cdot[H_L].
  \]
  Thus, we have proved $L$ splits $\xi$ under the additional
  hypothesis that $a\in L^{\times2}$.

  For the rest of the proof of~(i), assume $a\notin L^{\times2}$. Let
  $L'=L(\sqrt a) = L\otimes_kK$, and write again $s\colon L'\to L$ for
  the $L$-linear extension of $s$ to $L'$ and $N\colon L'\to L$ for
  the norm map. If $t\in L^{\times2}$, then
  $\xi_L=e_3(s_\star(\psi))_L$. Moreover, $L$ splits $H$ because it
  splits $U$. Therefore, \eqref{eq:e3psi} shows that $L$ splits
  $e_3(s_\star(\psi))$. For the rest of the proof, we may thus also
  assume $t\notin L^{\times2}$.

  Since $(a,t)_L=0$, we may find $z_0\in L'$ such that
  $t=N(z_0)$. Because $L$ splits $(b,t)_F$, we have $(b,N(z_0))_L=0$,
  so $\cores_{L'/L}(b,z_0)_{L'}=0$. It follows that $(b,z_0)_{L'}$ has
  an involution of the second kind, hence also a descent to $L$ by a
  theorem of Albert (see~\cite[(2.22)]{KMRT}). We may choose a descent
  of the form
  $(b,z_0)_{L'}=(b,\zeta)_{L'}$ for some $\zeta\in L^\times$; see
  \cite[(2.6)]{T}. Let
  $z=z_0\zeta\in {L'}^\times$. We then have $(b,z)_{L'}=0$, hence
  after taking the corestriction to $L$
  \[
  (b,N(z))_L=0.
  \]
  We also have $t=N(z_0)\equiv N(z)\bmod L^{\times2}$. Since $L$
  splits $[H_F]+(c,t)_F$, we have
  \[
  H_L=(c,N(z))_L.
  \]
  Since $H=(bc,N(x))_k=(c,N(y))_k$ by~(\ref{quat.eq}), it follows that
  \[
  (bc,N(xz))_L=(c,N(yz))_L=0.
  \]
  If $s(xz)=0$ (i.e., $xz\in L$), then
  $s_\star(\qform{xz}\pform{bc})$ is hyperbolic. If
  $s(xz)\neq0$, 
  computation yields 
  \[
  s_\star(\qform{xz}\pform{bc}) = \qform{s(xz)}\pform{bc, N(xz)};
  \]
  but since the quaternion algebra $(bc,N(xz))_{L}$ is split, the
  form $s_\star(\qform{xz}\pform{bc})$ is also hyperbolic in this
  case. Therefore, we may find $\lambda\in L^\times$ represented by
  $\qform{xz}\pform{bc}$; we then have
  \begin{equation}
    \label{eq:xi1}
    \qform{xz}\pform{bc} = \qform\lambda\pform{bc},\qquad\text{hence
      also}\quad \qform x\pform{bc}=\qform{\lambda z}\pform{bc}.
  \end{equation}
  Similarly, since the quaternion algebra $(c,N(yz))_{L}$ is split,
  the form $s_\star(\qform{yz}\pform c)$ is hyperbolic, and we may
  find $\mu\in L^\times$ such that
  \begin{equation}
    \label{eq:xi2}
    \qform{yz}\pform c = \qform\mu\pform c, \qquad\text{hence
      also}\quad \qform y\pform c=\qform{\mu z}\pform c.
  \end{equation}
  As a result of \eqref{eq:xi1} and \eqref{eq:xi2}, we have
  $\qform{x,-bcx} = \qform{\lambda z,-\lambda zbc}$ and $\qform{-y,cy}
  = \qform{-\mu z,\mu zc}$, hence we may rewrite $\psi$ over
  $L'$ as
  \[
  \psi_{L'}=\qform{bc, -c, \lambda z, -\lambda zbc, -\mu z,
    \mu zc}.
  \]
  Note that $z\notin L$ since $t\notin L^{\times2}$, hence
  $s(z)\neq0$. Using the last expression for $\psi_{L'}$ we may now compute 
  \[s_\star(\psi)_{L}  =s_\star(\psi_{L'})= s_\star(\qform z) \qform{\lambda,-\lambda bc,
    -\mu, \mu c}
   = \qform{s(z)}\pform{N(z)}\qform{\lambda,-\lambda bc,
    -\mu, \mu c}.
  \]
  Since $(bc,N(z))_{L}=(c,N(z))_{L}=H_{L}$, we have
  $\pform{N(z),bc} = \pform{N(z),c} = (n_H)_{L}$, hence
  $s_\star(\psi)_{L}=\qform{s(z)}\qform{\lambda,-\mu}(n_H)_{L}$,
  and therefore
  \begin{equation}
    \label{eq:xi3}
    e_3(s_\star(\psi))_{L} = (\lambda\mu)\cdot[H_{L}].
  \end{equation}
  On the other hand, we have $[C_K]=(bc,x)_K + (c,y)_K$, hence since
  $(b,z)_{L'}=0$
  \[
  [C_{L'}]= (bc,xz)_{L'} + (c,yz)_{L'}.
  \]
  In view of \eqref{eq:xi1} and \eqref{eq:xi2}, we may rewrite the
  right side as follows:
  \[
  [C_{L'}] = (bc,\lambda)_{L'} +
  (c,\mu)_{L'}.
  \]
  Therefore, $[C_{L}]+(bc,\lambda)_{L}+(c,y)_{L}$ is split by
  $L'$. We may then find $\nu\in L^\times$ such that
  \[
  [C_{L}] = (bc,\lambda)_{L} + (c,\mu)_{L} + (a,\nu)_{L}.
  \]
  Since $L$ splits $U$, we have $(t,a)_{L}=(t,b)_{L}=0$ and $(t,c)_{L} = H_{L}$. It 
  follows that
  \[
  (t)\cdot[C_{L}] = (t,c,\lambda\mu)_{L} =
  (\lambda\mu)\cdot[H_{L}].
  \]
  By comparing with \eqref{eq:xi3}, we see that $\xi$ vanishes over
  $L$. The proof of the lemma is thus complete. 
\end{proof}

\begin{proof}[Proof of Theorem~\ref{thm:xi}]
  Since $2\xi=0$, the assertion $f_3(U)=0$ follows from
  $\ch_U=\{0,\overline{\xi}^U\}$, see Remark~\ref{order.rem}(ii).   
 Moreover, the field $F_U$ splits $U$. Therefore, by Lemma~\ref{lem:tech}, we have $\xi\in\ker\bigl(H^3(F)\ra H^3(F_U)\bigr)$, so that 
 $\overline{\xi}^U\in\ch_U$. Since we know from
 Theorem~\ref{Peyre-cns} that the order of $\ch_U$ is at most $2$, it
 suffices to show that $\overline{\xi}^U\not =0$ when $\ch_U\not=0$ to
 establish $\ch_U=\{0,\overline{\xi}^U\}$. Therefore, proving the
 equivalence of (a), (b), (c) and (d) completes the proof.  

Let us first prove (a)~$\Rightarrow$~(b). Suppose $[C]=(a,\alpha)_k+(b,\beta)_k+(c,\gamma)_k$ for some
  $\alpha$, $\beta$, $\gamma\in k^\times$. Since
  $[C_K]=(bc,x)_K+(c,y)_K$, it follows that $(b,\beta)_K+(c,\gamma)_K=
  (bc,x)_K+(c,y)_K$, hence
  \[
  (bc,\beta x)_K=(c,\beta\gamma y)_K.
  \]
  By the common slot lemma \cite[Th.~III.4.13]{Lam}, we may find
  $z\in K^\times$ such that
  \begin{equation}
    \label{eq:tech1}
    (bc,\beta x)_K=(bc,z)_K=(c,z)_K=(c,\beta\gamma y)_K.
  \end{equation}
  Let $E=F\bigl(\sqrt{N(z)t}\,\bigr)$, a quadratic extension of $F$. We claim that
  $E$ splits $U$. First, observe that $N(z)t$ is represented by the
  form $\qform{t,-at}$, hence the quaternion algebra $(a,t)_F$
  contains a pure quaternion with square $N(z)t$. Therefore, $E$
  splits $(a,t)_F$. Likewise, from~\eqref{eq:tech1} we see that
  $(b,z)_K=0$, hence by taking the corestriction to $k$ we have
  $(b,N(z))_k=0$. Therefore, $N(z)t$ is represented by the form
  $\qform{t,-bt}$, and it follows that $E$ splits $(b,t)_F$. Finally,
  by taking the corestriction of each side of the rightmost equation
  in~\eqref{eq:tech1}, we obtain $(c,N(z))_k=(c,N(y))_k$, so
  $N(y)N(z)$ is represented by $\qform{1,-c}$ and therefore $N(z)t$ is
  represented by $\qform{N(y)t,-cN(y)t}$. It follows that $E$ splits
  the quaternion algebra $(c,N(y)t)_F$. We have thus shown that $E$
  splits three generators of $U$, hence $E$ splits $U$.

  The implication (b)~$\Rightarrow$~(c) follows immediately from
  Peyre's Theorem~\ref{Peyre-cns}.  Moreover, 
  (c)~$\Rightarrow$~(d) is clear since $\overline{\xi}^U\in\ch_U$. 
  To complete the proof, we show (d)~$\Rightarrow$~(a). Suppose there
  exist $\lambda_1$, $\lambda_2$, $\lambda_3\in F^\times$ such that
  \begin{equation}
    \label{eq:xi10}
  \xi=(\lambda_1,a,t)+(\lambda_2,b,t)+(\lambda_3,c,N(y)t).
  \end{equation}
  Let $\partial\colon H^i(F)\to H^{i-1}(k)$ be the residue map
  associated to the $t$-adic valuation, for $i=2$, $3$. Since
  $e_3(s_\star(\psi))\in H^3(k)$ we have
  $\partial(e_3(s_\star(\psi)))=0$, hence
  $\partial(\xi)=[C]$. Therefore, taking the image of each side
  of~\eqref{eq:xi10} under the residue map yields
  \[
  [C]=a\cdot\partial(\lambda_1,t) + b\cdot\partial(\lambda_2,t) +
  c\cdot\partial(\lambda_3,N(y)t),
  \]
  so that $[C]\in \Dec(M/k)$. 
\end{proof}

As a corollary, we get:

\begin{cor}
  \label{decsub.cor}
  Use the notation~\ref{not:not}, and assume
  $[C]\notin\Dec(M/k)$. Then $U\subset\br(F)$ is a quaternionic
  subgroup of order~$8$ such that $\sum_{[H]\in U} n_H\in I^4(F)$
  (i.e., $f_3(U)=0$), which is not split by any quadratic extension
  of $F$ (i.e., $\ch_{U}\neq0$).
\end{cor}

To obtain an example of a central simple algebra with orthogonal
involution of degree~$12$ without quadratic splitting, we need a more
stringent condition on $C$:

\begin{lem}
  \label{lem:stringent}
  With the notation~\ref{not:not}, the following conditions are
  equivalent:
  \begin{enumerate}
  \item[(a)]
  there is a quadratic extension $E$ of $F$ that splits $(a,t)_F$ and
  $\xi$;
  \item[(b)]
  the algebra $C$ is Brauer-equivalent to a tensor product of
  quaternion $k$-algebras $A_1\otimes_kA_2\otimes_kA_3$ with $A_3$
  split by $K$.
  \end{enumerate}
\end{lem}

\begin{proof}
  (a)~$\Rightarrow$~(b):
  Let $\widehat F=k((t))$ be the completion of $F$ for the
  $t$-adic valuation. The field $E$ does not embed in $\widehat F$
  because $\widehat F$ does not split $(a,t)_F$. Therefore, $E$ and
  $\widehat F$ are linearly disjoint over $F$ and we may consider the
  field $\widehat E=E\otimes_F\widehat F$, which is a quadratic
  extension of $\widehat F$ that splits $(a,t)_{\widehat F}$ and
  $\xi_{\widehat F}$. Each square class in $\widehat F$ is represented
  by an element in $k^\times$ or an element of the form $ut$ with
  $u\in k^\times$, see \cite[Cor.~VI.1.3]{Lam}. Therefore, we may
  assume that either 
  $\widehat E=\widehat F(\sqrt{u})$ or $\widehat E=\widehat
  F(\sqrt{ut})$ for some $u\in k^\times$.

  Suppose first $\widehat E=\widehat F(\sqrt u)$ with $u\in
  k^\times$. Since the quaternion algebra $(a,t)_{\widehat F}$ is
  split by $\widehat E$, it must contain a pure quaternion with square
  $u$, hence $u$ is represented by $\qform{a,t,-at}$ over $\widehat
  F$. Therefore, $u\equiv a\bmod k^{\times2}$, and $\widehat
  E=K((t))$. From
  \[
  \xi_{\widehat E} = t\cdot[C_{\widehat E}] +
  e_3(s_\star(\psi))_{\widehat E}=0,
  \]
  it follows by taking images under the residue map $H^3(\widehat
  E)\to H^2(K)$ associated to the $t$-adic valuation that
  $[C_K]=0$. Then $C$ is Brauer-equivalent to a quaternion algebra
  $A_3$ split by $K$, and (b) holds with $A_1$, $A_2$ split quaternion
  algebras. 

  Suppose next $\widehat E=\widehat F(\sqrt{ut})$ for some $u\in
  k^\times$. Since $\widehat E$ splits $(a,t)_{\widehat F}$, it
  follows as above that $ut$ is represented by $\qform{a,t,-at}$ over
  $\widehat F$, hence $u$ is represented by $\qform{1,-a}$, which
  means that $u\in N(K^\times)$. Because $ut$ is a square in $\widehat
  E$, we have $t\cdot[C_{\widehat E}]=u\cdot[C_{\widehat E}]$, hence
  the equation $\xi_{\widehat E}=0$ yields
  \[
  u\cdot[C_{\widehat E}] + e_3(s_\star(\psi))_{\widehat E} =
  \bigl(u\cdot[C]+e_3(s_\star(\psi))\bigr)_{\widehat E}=0.
  \]
  Since $\widehat F=k((t))=k((ut))$ we have $\widehat
  E=k((\sqrt{ut}))$, hence the scalar extension map $H^3(k)\to
  H^3(\widehat E)$ is injective. Therefore, the last equation yields
  \[
  u\cdot[C]+e_3(s_\star(\psi))=0,
  \]
  which shows that $e_3(s_\star(\psi))\in N(K^\times)\cdot[C]$ because
  $u\in N(K^\times)$. By Barry's result \cite[Prop.~4.4]{Barry}, it
  follows that the 
  biquaternion algebra $B$ has a descent to $k$: there exist
  quaternion $k$-algebras $A_1$, $A_2$ such that $B\simeq
  A_1\otimes_kA_2\otimes_kK$. 
  Since $C_K$ is Brauer-equivalent to $B$, it follows that $C\otimes
  A_1\otimes_kA_2$ is split by $K$. It is therefore
  Brauer-equivalent to a quaternion algebra $A_3$ split by $K$, and
  $C$ is Brauer-equivalent to $A_1\otimes_kA_2\otimes_kA_3$, proving~(b).
  \smallbreak  

  (b)~$\Rightarrow$~(a):
  Since $B$ is Brauer-equivalent to $C_K$, condition~(b) implies that
  $B\simeq A_1\otimes_kA_2\otimes_kK$. From Barry's result
  \cite[Prop.~4.4]{Barry}, it follows that $e_3(s_\star(\psi))=
  u\cdot[C]$ for some $u\in N(K^\times)$. Let $E=F(\sqrt{ut})$. Then
  $(a,t)_E\simeq (a,u)_E$, hence $E$ splits $(a,t)_F$ because $u\in
  N(K^\times)$. Moreover,
  $\xi_E=\bigl(u\cdot[C]+e_3(s_\star(\psi))\bigr)_E$, so $E$ also
  splits $\xi$. Therefore, (a) holds.
\end{proof}

Examples of algebras $C$ for which condition~(b)
of Lemma~\ref{lem:stringent} does not hold include indecomposable
division algebras of degree~$8$ and exponent~$2$; other examples are
given in \cite{Barry2}. Note that condition~(b) is weaker than
$[C]\in\Dec(M/k)$; it is in fact strictly weaker: see
Remark~\ref{rem:weaker}. 

\begin{cor}
  \label{cor:comeback}
  Use the notation~\ref{not:not}, and let
  $Q=(a,t)_F$. There exists an orthogonal involution $\rho$ on
  $M_6(Q)$ with the following properties:
  \begin{enumerate}
  \item[(i)] $\rho$ has trivial discriminant and trivial Clifford
    invariant;
  \item[(ii)] $e_3(\rho)=\xi\mod F^\times\cdot[Q]$;
  \item[(iii)] $U$ is a decomposition group of $\rho$;
  \item[(iv)] $f_3(\rho)=0$.
  \end{enumerate}
  For any involution $\rho$ satisfying~\emph{(i)} and \emph{(ii)},
  there exists a 
  quadratic extension 
  of $F$ over which $Q$ is split and $\rho$ is hyperbolic if and only
  if the equivalent conditions of Lemma~\ref{lem:stringent} hold. 
\end{cor}

\begin{proof}
  By Remark~\ref{rem:anyU}, there is an orthogonal involution $\rho$
  on $M_6(Q)$ with trivial discriminant and trivial Clifford
  invariant, and with decomposition group~$U$. By
  Theorems~\ref{hom.thm} and~\ref{f3.thm}, $\overline{e_3(\rho)}^U$
  generates $\ch_U$, 
  and $f_3(\rho)=f_3(U)$. Therefore, Theorem~\ref{thm:xi} yields
  $\overline{e_3(\rho)}^U=\overline{\xi}^U$ and $f_3(\rho)=0$. By
  Proposition~\ref{same-dec.prop}(ii), we may assume
  $e_3(\rho)=\xi\mod F^\times\cdot[Q]$. Thus, $\rho$ satisfies
  conditions~(i)--(iv). 

  Now, let $\rho$ be any orthogonal involution on $M_6(Q)$
  satisfying~(i) and (ii). Because of~(ii), condition~(a) of
  Lemma~\ref{lem:stringent} holds if and only if there is a 
  quadratic extension $E$ of $F$ such that $[Q_E]=0$ and
  $e_3(\rho)_E=0$. By Theorem~\ref{iso.thm}(i), the last equation
  holds if and only if $\rho_E$ is hyperbolic.
\end{proof}

Suppose $C$ does not satisfy condition~(b) of
Lemma~\ref{lem:stringent} (e.g., $C$ is an indecomposable division
algebra of degree~$8$ and exponent~$2$ split by $M$). Then for any
involution $\rho$ on $M_6(Q)$ satisfying the properties~(i)--(iv) of
Corollary~\ref{cor:comeback} we have $f_3(\rho)=0$, and yet there is
no quadratic extension of $F$ over which $Q$ is split and $\rho$ is
hyperbolic. From (b')~$\Rightarrow$~(a') in
Proposition~\ref{f3zero.prop}, it follows that $\ch_{U'}\neq0$ for every
decomposition group $U'$ of $\rho$.

\begin{rem}
  \label{rem:weaker}
  By \cite[Cor.~3.2]{T}, for any triquadratic extension $M/k$, any
  $2$-torsion Brauer class in $\Br(k)$ split by $M$ is represented
  modulo~$\Dec(M/k)$ by a quaternion algebra. Therefore, if the
  triquadratic extension $M/k$ is such that $\Dec(M/k)$ does not
  coincide with the subgroup of ${}_2\Br(k)$ split by $M$ (see
  \cite[\S5]{ELTW} for examples of such extensions), we may find a
  quaternion $k$-algebra $C$ split by $M$ such that
  $[C]\notin\Dec(M/k)$. The algebra $C$ obviously satisfies
  condition~(b) of Lemma~\ref{lem:stringent} (with $A_2$ and $A_3$
  split), so for any involution $\rho$ on $M_6(Q)$ satisfying the
  properties~(i)--(iv) of Corollary~\ref{cor:comeback} we may find a
  quadratic extension of $F$ over which $Q$ is split and $\rho$ is
  hyperbolic. From~(a')~$\Rightarrow$~(b') in
  Proposition~\ref{f3zero.prop}, it follows that there exists a
  decomposition group $U'$ of $\rho$ such that $\ch_{U'}=0$. Yet,
  because $[C]\notin\Dec(M/k)$, the decomposition group $U$ of $\rho$
  satisfies $\ch_U\neq0$ by Theorem~\ref{thm:xi}.
\end{rem}

\section{Application to degree $8$ algebras with involution}

The Arason invariant in degree $8$ was studied in~\cite{QT:Arason} for orthogonal involutions with trivial discriminant and trivial Clifford algebra. 
In this section, we extend it to algebras of degree $8$ and index $2$,
when the involution has trivial discriminant and the two components of
the Clifford algebra also both have index $2$.  
First, we prove an analogue of Theorem~\ref{thm:adddec12} on additive
decompositions, for degree $8$ algebras with orthogonal involution of
trivial discriminant. 

\subsection{Additive decompositions in degree $8$}

\label{sec:adddec8}

Let $(A,\sigma)$ be a degree $8$ algebra with orthogonal involution of trivial discriminant. 
We let $(C^+(A,\sigma),\sigma^+)$ and
$(C^-(A,\sigma),\sigma^-)$ denote the two components of the Clifford
algebra of $(A,\sigma)$, endowed with the involutions induced by the
canonical involution of the Clifford algebra.  
Both algebras have degree $8$, both involutions have trivial discriminant, and by triality~\cite[(42.3)]{KMRT}, 
\begin{equation}
  \label{eq:Ctri}
  C(C^+(A,\sigma),\sigma^+)\simeq (C^-(A,\sigma),\sigma^-)\times
  (A,\sigma)
  \end{equation}
  and
  \begin{equation}
  \label{eq:Ctri2}
  C(C^-(A,\sigma),\sigma^-) \simeq
  (A,\sigma) \times (C^+(A,\sigma),\sigma^+).
  \end{equation}

Assume $(A,\sigma)$ decomposes as a sum $(A,\sigma)\in
  (A_1,\sigma_1)\orthsum(A_2,\sigma_2)$ of two degree $4$ algebras with orthogonal involution of trivial discriminant. Each summand is a tensor product of two quaternion algebras with canonical involution, and we get 
 \begin{equation}
 \label{eq:dec}
 (A,\sigma)\in\bigl((Q_1,\ba)\otimes(Q_2,\ba)\bigr)\boxplus\bigl((Q_3,\ba)\otimes(Q_4,\ba)\bigr),
 \end{equation}
  for some quaternion algebras $Q_1$, $Q_2$, $Q_3$ and $Q_4$ such that
  $A$ is Brauer-equivalent to $Q_1\otimes Q_2$ and $Q_3\otimes Q_4$.  
 By~\cite[Prop. 6.6]{QSZ}, the two components of the Clifford algebra of $(A,\sigma)$ then admit similar decompositions, namely, up to permutation of the two components, we have: 
\begin{equation} 
\label{eq:c+}
\bigl(C^+(A,\sigma),\sigma^+\bigr)\in\bigl((Q_1,\ba)\otimes(Q_3,\ba)\bigr)\boxplus\bigl((Q_2,\ba)\otimes(Q_4,\ba)\bigr),
\end{equation} and 
\begin{equation}
\label{eq:c-}
\bigl(C^-(A,\sigma),\sigma^-\bigr)\in\bigl((Q_1,\ba)\otimes(Q_4,\ba)\bigr)\boxplus\bigl((Q_2,\ba)\otimes(Q_3,\ba)\bigr).
\end{equation}

Mimicking the construction in \S\ref{sec:adddec12}, we associate to
every decomposition of $(A,\sigma)$ as above the subgroup $W$ of the Brauer group of $F$ generated by any three elements among the $[Q_i]$ for $1\leq i\leq 4$. We call $W$ a decomposition group of $(A,\sigma)$. It consists of at most $8$ elements, and can be described explicitly by 
\[W=\{0,[A],[C^+(A,\sigma)],[C^-(A,\sigma)],[Q_1],[Q_2],[Q_3],[Q_4]\}.\] 
In view of their additive decompositions, $W$ also is a decomposition group of the two components 
$(C^+(A,\sigma),\sigma^+)$ and
$(C^-(A,\sigma),\sigma^-)$ of the Clifford algebra. 
Note that, in contrast with the decomposition groups of algebras of
degree~$12$ in Definition~\ref{defi:decsubgrp}, the group $W$ may
contain three Brauer classes of index~$4$ instead of at most
one. Nevertheless, it has similar properties; for instance, we prove: 

\begin{prop}
  \label{prop:quadsplit8}
  Let $(A,\sigma)$ be a central simple algebra of degree~$8$ with an
  orthogonal involution of trivial discriminant. 

\begin{enumerate}
\item[(i)] Suppose $(A,\sigma)$ has an additive decomposition as
  in~\eqref{eq:dec}, with decomposition group $W$.  For any extension
  $K/F$ which splits $W$, the algebra with involution $(A_K,\sigma_K)$
  is split and hyperbolic.   
\item[(ii)] The converse holds for quadratic extensions: if $(A,\sigma)$ is
    split and hyperbolic over a quadratic extension $K$ of $F$, then
    $(A,\sigma)$ has an additive decomposition with decomposition
    group split by $K$.
  \end{enumerate}

 \end{prop}

\begin{proof}
(i)
  If a field $K$ splits $W$,
  then it splits $A$, and moreover each summand in~(\ref{eq:dec}) is
  split and hyperbolic over $K$, therefore  
  $\sigma_K$ is hyperbolic.

(ii)
  To prove the converse, suppose $K=F(d)$ with $d^2=\delta\in
  F^\times$, and assume $A_K$ is split and $\sigma_K$ is hyperbolic,
  hence $\ind A\leq2$. If $A$ is split, we have as in the proof of
  Proposition~\ref{f3zero.prop} $(A,\sigma)\simeq\Ad_\varphi$ with
  $\varphi$ an $8$-dimensional quadratic form multiple of
  $\qform{1,-\delta}$. We may then find quaternion $F$-algebras $Q_1$,
  $Q_2$ split by $K$ and scalars $\alpha_1$, $\alpha_2\in F^\times$
  such that $\varphi\simeq\qform{\alpha_1}n_{Q_1}\perp
  \qform{\alpha_2}n_{Q_2}$. As in Example~\ref{split.ex}, we obtain a
  decomposition 
  \[
  (A,\sigma)\in \bigl((Q_1,\ba)\otimes(Q_1,\ba)\bigr) \orthsum
  \bigl((Q_2,\ba)\otimes(Q_2,\ba)\bigr).
  \]
  The corresponding decomposition group is
  $\{0,[Q_1],[Q_2],[Q_1]+[Q_2]\}$; it is split by $K$.

  If $\ind A=2$, let $Q$ be the quaternion division algebra
  Brauer-equivalent to $A$. Since $K$ splits $A$, we may, again as in the
  proof of Proposition~\ref{f3zero.prop}, identify $K$ with a
  subfield of $Q$ and find a skew-hermitian form $h$ of the form
  $\qform d\qform{1,\alpha,\beta,\gamma}$ (with $\alpha$, $\beta$,
  $\gamma\in F^\times$) such that $(A,\sigma)\simeq\Ad_h$. Then
  \[
  (A,\sigma)\in\Ad_{\qform d\qform{1,\alpha}} \orthsum \Ad_{\qform
    d\qform{\beta,\gamma}}
  \]
  is a decomposition in which each of the summands becomes hyperbolic
  over $K$. The corresponding decomposition group is therefore split
  by $K$.
\end{proof}

There exist quadratic forms $\varphi$ of dimension~$8$ with trivial
  discriminant and Clifford algebra of index~$4$ that do not decompose
  into an orthogonal sum of two $4$-dimensional quadratic forms of
  trivial discriminant, see \cite[Cor.~16.8]{IzhKar} or
  \cite[Cor~6.2]{HofTig}. For such a form, neither $\Ad_\varphi$ nor the components of its Clifford algebra have additive decompositions as in~(\ref{eq:dec}). 
  The next proposition shows, by contrast, that such a decomposition always exist if at least two among the algebras $A$, $C^+(A,\sigma)$ and $C^-(A,\sigma)$ have index $\leq 2$.   
    
  \begin{prop} 
  \label{prop:adddec8} 
  Let $(A,\sigma)$ be a central simple algebra of degree $8$ with orthogonal involution of trivial discriminant. We assume at least two among the algebras $A$, $C^+(A,\sigma)$ and $C^-(A,\sigma)$ have index $\leq 2$. Then all three algebras with involution $(A,\sigma)$, $(C^+(A,\sigma),\sigma^+)$ and $(C^-(A,\sigma),\sigma^-)$ admit an additive decomposition as a sum of two degree $4$ algebras with orthogonal involution of trivial discriminant as in~(\ref{eq:dec}).   
  \end{prop} 
  
  \begin{proof}
  
  Assume two among $\ind A$, $\ind C^+(A,\sigma)$, $\ind
  C^-(A,\sigma)$ are $1$ or $2$. By triality, see (12) to (16) above, it is enough to prove that one of the three algebras with involution, say $(A,\sigma)$ has an additive decomposition. 
  Since $A$, $C^+(A,\sigma)$,
  $C^-(A,\sigma)$ are interchanged by triality, we may also assume $\ind
  A\leq2$. If $A$ is split, so $(A,\sigma)\simeq\Ad_\varphi$ for some
  $8$-dimensional quadratic form $\varphi$ with trivial discriminant
  and Clifford algebra of index at most~$2$, then (a) holds by a
  result of Knebusch \cite[Ex.~9.12]{Knebusch}, which shows that
  $\varphi$ is the product of 
  a $2$-dimensional quadratic form and a $4$-dimensional quadratic
  form. 

  For the rest of the proof, assume $(A,\sigma)\simeq\Ad_h$
  for some skew-hermitian form $h$ over a quaternion division algebra
  $(Q,\ba)$. Let $q\in Q$ be a nonzero quaternion represented by $h$,
  and let $h\simeq\qform q \perp h'$ for some skew-hermitian form $h'$
  of absolute rank~$6$. As we saw in the proof of
  Theorem~\ref{thm:adddec12}, over the quadratic extension $K=F(q)$
  the algebra $Q$ splits and the form $\qform q$ becomes hyperbolic,
  hence $h_K$ and $h'_K$ are Witt-equivalent. In particular, it
  follows that $e_2((\ad_{h'})_K)= e_2((\ad_h)_K)$ has index at
  most~$2$. But $(\ad_{h'})_K=\ad_\psi$ for some Albert form~$\psi$
  over $K$, so $\psi$ is isotropic. It follows by \cite[Prop.,
  p.~382]{QT} that $h'$ 
  represents some scalar multiple of $q$; thus $h\simeq \qform q
  \qform{1,-\lambda}\perp h''$ for some $\lambda\in F^\times$ and some
  skew-hermitian form $h''$ of absolute rank~$4$. The discriminant of
  $h''$ must be trivial because $h$ and $\qform q\qform{1,-\lambda}$
  have trivial discriminant, and we thus have the required
  decomposition for $(A,\sigma)$.  
\end{proof}

\begin{rem}
It follows that all trialitarian triples such that at least two of the algebras have index $\leq 2$ have a description as 
in~(\ref{eq:dec}) to (\ref{eq:c-}). 
\end{rem} 
\subsection{An extension of the Arason invariant in degree $8$ and index
  $2$.}

Throughout this section, $(A,\sigma)$ is a central simple $F$-algebra
of degree~$8$ and trivial discriminant. It is known that $(A,\sigma)$
is a tensor product of quaternion algebras with involution if and only
if $e_2(\sigma)=0$, see \cite[(42.11)]{KMRT}. In this case, the Arason
invariant 
$e_3(\sigma)\in M^3_A(F)$ is defined when $A$ has index at most $4$ (see \S\ref{subsec:Arasoninvol}) , and
represented by an element of order~$2$ in $H^3(F)$, see
\cite{QT:Arason}. Here, we
extend the definition of the $e_3$ invariant under the following
hypothesis: 
\begin{equation}
\label{eq:C}
\ind A=\ind C^+(A,\sigma) = \ind C^-(A,\sigma) = 2.
\end{equation}
By Proposition~\ref{prop:adddec8}, this condition implies that
$(A,\sigma)$ decomposes into a sum of two central simple algebras of
degree~$4$ with involutions of trivial discriminant. Moreover, the
associated decomposition group $W$ is a quaternionic subgroup of
$\br(F)$. 
Let $Q$, $Q^+$, $Q^-$ be the quaternion division algebras over $F$
that are Brauer-equivalent to $A$, $C^+(A,\sigma)$, and
$C^-(A,\sigma)$ respectively. 
>From the Clifford algebra relations
\cite[(9.12)]{KMRT}, 
we know $[Q^+]+[Q^-]=[Q]$. Therefore, the following is a subgroup of
the Brauer group:
\[
V=\{0,[Q], [Q^+], [Q^-]\}\subset\br(F).
\]
Condition~\eqref{eq:C} implies that $\vert V\rvert=4$. Moreover, $V$
also is a subgroup of every decomposition group of $(A,\sigma)$.

To $(A,\sigma)$, we may associate algebras of degree $12$ with orthogonal involution with trivial discriminant and trivial Clifford invariant by considering any involution $\rho$ of $M_6(Q)$ such that 
 \begin{equation}
    \label{eq:rho}
    (M_6(Q),\rho)\in (A,\sigma)\orthsum\bigl((Q^+,\ba)\otimes
    (Q^-,\ba)\bigr). 
  \end{equation}
  Since the two components of the Clifford algebra of $\sigma$ are Brauer-equivalent to $Q^+$ and $Q^-$, 
  the involution $\rho$ has trivial Clifford invariant. Therefore, we may consider its Arason invariant 
  $e_3(\rho)\in M_Q^3(F)$. 
 The following lemma compares the Arason invariant of two such involutions: 
 \begin{lem}
 Let $\rho$ and $\rho'$ be two involutions of $M_6(Q)$ satisfying~(\ref{eq:rho}). 
 There exists $\lambda\in F^\times$ such that 
 \[e_3(\rho)-e_3(\rho')=(\lambda)\cdot [Q^+]=(\lambda)\cdot[Q^-]\in M_Q^3(F).\]
 Moreover, $f_3(\rho)=f_3(\rho')$. 
 \end{lem} 
  \begin{proof} 
  
By definition of the direct orthogonal sum for algebras with involution, 
we may pick skew-hermitian forms $h_1$ and $h_2$ over $(Q,\ba)$ such that $\sigma=\ad_{h_1}$, 
$\ba\otimes\ba\simeq \ad_{h_2}$ and $\rho=\ad_{h_1\perp h_2}$. 
Moreover, there exists $\lambda\in F^\times$ such that $\rho'=\ad_{h_1\perp\qform{\lambda}h_2}$. 
Therefore, we have (see~\S\,\ref{subsec:relAra} and ~\ref{subsec:Arasoninvol}): 
\[e_3(\rho)-e_3(\rho')=e_3(\qform{1,-\lambda}h_2)=(\lambda)\cdot [Q^+]\mod F^\times\cdot [Q],\]
by proposition~\ref{compute.prop}, since $e_2(h_2)=[Q^+]=[Q^-]\mod[Q]$. 
Moreover, if $c$ and $c'\in H^3(F)$ are representatives of $e_3(\rho)$ and $e_3(\rho')$ respectively, then 
\[c'-c\in (\lambda)\cdot[Q^+]+F^\times \cdot [Q].\] So $2c=2c'$, and this finishes the proof. 

\end{proof} 

Let $F^\times\cdot V\subset H^3(F)$ be the subgroup
consisting of the products $\lambda\cdot v$ with $\lambda\in F^\times$
and $v\in V=\{0,[Q],[Q^+],[Q^-]\}$. This subgroup contains $F^\times\cdot[Q]$, so we
may consider the canonical map $\ba^V\colon M^3_Q(F)\to
H^3(F)/F^\times\cdot V$. The previous lemma shows that the image $\overline{e_3(\rho)}^V$ of the Arason invariant of $\rho$ does not depend on the choice of an involution $\rho$ satisfying~(\ref{eq:rho}). 
This leads to the following: 

\begin{defi}
  \label{defi:e3deg8}
  With the notation above, we set
  \[
  e_3(\sigma) = \overline{e_3(\rho)}^V\in H^3(F)/F^\times\cdot V
  \qquad\text{and}\qquad 
  f_3(\sigma) = f_3(\rho)\in F^\times\cdot[Q]\subset H^3(F)
  \]
 where $\rho$ is any involution satisfying
  \begin{equation*}
       (M_6(Q),\rho)\in (A,\sigma)\orthsum\bigl((Q^+,\ba)\otimes
    (Q^-,\ba)\bigr). 
  \end{equation*}
    \end{defi}

This definition functorially extends the definition of the Arason invariant. 
Indeed, 
if $K$ is any extension of $F$ that splits $Q^+$ or $Q^-$ (or both),
then the scalar extension map $\br(F)\to\br(K)$ carries $V$ to
$\{0,[Q]\}$, and any involution $\rho$ as in~\eqref{eq:rho} becomes
Witt-equivalent to $\sigma$ over $K$. Therefore, scalar extension
carries $e_3(\sigma)\in H^3(F)/F^\times\cdot V$ defined above to
$e_3(\sigma_K)\in M^3_Q(K)$ as defined in
\S\ref{subsec:Arasoninvol}. 

\begin{ex}
Consider a central simple algebra $(M_6(Q),\rho)$ of degree~$12$ and
index~$2$ with an orthogonal involution of trivial discriminant and
trivial Clifford invariant. By Theorem~\ref{thm:adddec12}
$(M_6(Q),\rho)$ admits additive decompositions 
\[
(M_6(Q),\rho)\in \orthsum_{i=1}^3\bigl((Q_i,\ba)\otimes(H_i,\ba)
\bigr)\qquad\text{with $\sum_{i=1}^3[H_i]=0$},
\]
so it contains symmetric idempotents $e_1$, $e_2$, $e_3$ such that
\[
\bigl(e_iM_6(Q)e_i, \rho\rvert_{e_iM_6(Q)e_i} \bigr) \simeq
(Q_i,\ba)\otimes(H_i,\ba).
\]
Consider the restriction of $\rho$ to
$(e_1+e_2)M_6(Q)(e_1+e_2)$; we thus obtain an algebra with involution
$(M_4(Q),\sigma)$ such that
\begin{multline*}
(M_6(Q),\rho)\in(M_4(Q),\sigma) \orthsum \bigl((Q_3,\ba)\otimes
(H_3,\ba)\bigr)\\ \quad\text{and}\quad
(M_4(Q),\sigma)\in \orthsum_{i=1}^2\bigl((Q_i,\ba)\otimes(H_i,\ba)
\bigr).
\end{multline*}
It is clear that the discriminant of $\sigma$ is trivial. Since
$e_2(\rho)=0$, we have 
\[
e_2(\sigma)=e_2\bigl((Q_3,\ba)\otimes(H_3,\ba)
\bigr)=\{[Q_3], [H_3]\}.
\]
Therefore, Condition~\eqref{eq:C} holds for $(M_4(Q),\sigma)$ if $Q_3$
and $H_3$ are not split. In that case, we have
$V=\{0,[Q],[Q_3],[H_3]\}$ and, by definition,
\[
e_3(\sigma)=\overline{e_3(\rho)}^V\in H^3(F)/F^\times\cdot V
\quad\text{and}\quad
f_3(\sigma)=f_3(\rho)\in F^\times\cdot[Q]\subset H^3(F).
\]
The condition that $Q_3$ and $H_3$ are not split holds in particular
when the decomposition group $U$ generated by $[Q_1]$, $[Q_2]$,
$[Q_3]$ has order~$8$.
\end{ex}
\begin{ex}
Take for $(M_6(Q),\rho)$ the algebra
$\Ad_{\qform{1,-t}}\otimes(\lambda^2E,\gamma)_{F(t)}$ of
Corollary~\ref{f3.cor}, with $E$ a division algebra of degree and
exponent~$4$. (Note that $\lambda^2E$ is Brauer-equivalent to
$E\otimes E$, hence it has index~$2$.) Since $f_3(\rho)\neq0$, every
decomposition group of $\rho$ has order~$8$; indeed, quaternionic
subgroups $U\subset\br(F)$ of order dividing $4$ have $\ch_U=0$ by
Theorem~\ref{Peyre}, hence trivial $f_3$ by
Proposition~\ref{hufu.prop}.  
The construction in the previous example yields an algebra
with involution $(M_4(Q),\sigma)$ of degree~$8$ satisfying
Condition~\eqref{eq:C}, with $f_3(\sigma)=t\cdot[Q]\neq0$.
\end{ex}
\begin{ex}
Also, we may take for $(M_6(Q),\rho)$ the algebra with involution of
Corollary~\ref{cor:comeback}, and obtain an algebra with involution
$(M_4(Q),\sigma)$ of degree~$8$ satisfying Condition~\eqref{eq:C} such
that (with the notation~\ref{not:not})
$e_3(\sigma)=\overline{\xi}^V\in H^3(F)/F^\times V$. Since $\xi\notin
F^\times\cdot U$, we have $e_3(\sigma)\neq0$. Yet, we have
$f_3(\sigma)=f_3(\rho)=0$ by Corollary~\ref{cor:comeback}. Moreover,
there is no quadratic extension $K$ of $F$ such that $Q_K$ is split
and $\sigma_K$ is hyperbolic. Indeed, over such a field, $(M_6(Q),\rho)_K$ would be
Witt-equivalent to an algebra of degree~$4$, hence it would be
hyperbolic because $e_2(\rho)=0$. Corollary~\ref{cor:comeback} shows
that such quadratic extensions $K$ do not exist.
\end{ex}

The next proposition shows that the $e_3$ invariant detects isotropy,
for any central simple algebra with involution $(A,\sigma)$ satisfying
Condition~\eqref{eq:C}. As in \S\ref{sec:adddec8}, we let $\sigma^+$
and $\sigma^-$ denote the canonical involutions on $C^+(A,\sigma)$ and
$C^-(A,\sigma)$. 

\begin{prop}
  \label{prop:12}
  Let $(A,\sigma)$ be a central simple $F$-algebra of degree~$8$ with
  orthogonal involution of trivial discriminant
  satisfying~\eqref{eq:C}. 
  With the notation above, we have
  $e_3(\sigma)=e_3(\sigma^+)=e_3(\sigma^-)$ and
  $f_3(\sigma)=f_3(\sigma^+) = f_3(\sigma^-)$. Moreover, the following
  conditions are equivalent:
  \begin{enumerate}
  \item[(a)]
  $e_3(\sigma)=0$;
  \item[(b)]
  $\sigma$ is isotropic;
  \item[(c)]
  $(A,\sigma)$ is Witt-equivalent to $(Q^+,\ba)
  \otimes (Q^-,\ba)$. 
  \end{enumerate}
\end{prop}

\begin{proof}
  As in \S\ref{sec:Peyre}, let $F_V$ denote the function field of the
  product of the Severi--Brauer varieties associated to the elements
  of $V$.
  Extending scalars to $F_V$, we split $Q$ and $e_2(\sigma)$, hence
  there is a $3$-fold Pfister form $\pi$ over $F_V$ such that
  \[
  \sigma_{F_V} \simeq \ad_\pi.
  \]
  For Pfister forms, we have $\ad_\pi\simeq\ad_\pi^+\simeq\ad_\pi^-$
  (see \cite[(42.2)]{KMRT}), hence
  $\sigma^+_{F_V} \simeq \sigma^-_{F_V} \simeq
  \sigma_{F_V}$, and therefore
  \[
  e_3(\sigma)_{F_V} = e_3(\sigma^+)_{F_V} = e_3(\sigma^-)_{F_V} =
  e_3(\pi).
  \]
  Since $V$ is generated by the Brauer classes of two quaternion
  algebras, it follows from Theorem~\ref{Peyre} that $F^\times\cdot V$
  is the kernel of the scalar extension map $H^3(F)\to H^3(F_V)$,
  hence the preceding equations yield
  $e_3(\sigma) = e_3(\sigma^+) = e_3(\sigma^-)$. We then have
  $f_3(\sigma)=f_3(\sigma^+)= f_3(\sigma^-)$, since $f_3(\sigma)$
  (resp.\ $f_3(\sigma^+)$, resp.\ $f_3(\sigma^-)$) is $2$ times any
  representative of $e_3(\sigma)$ (resp.\ $e_3(\sigma^+)$, resp.\
  $e_3(\sigma^-)$) in $H^3(F)$.

  To complete the proof, we show that (a), (b), and (c) are
  equivalent. 
  Clearly, (c) implies (b). The converse follows easily from~\cite[(15.12)]{KMRT} if $A$ has index $2$, and ~\cite[(16.5)]{KMRT} if $A$ is split. 
  Moreover, in view of the definition of $e_3(\sigma)$, the equivalence between (a) and (c) follows from Proposition~\ref{isotAra.prop}. 
 \end{proof}

As in \S\ref{sec:ArasonPeyre}, we may relate the $e_3$ invariant to
the homology of the Peyre complex of any decomposition group, as
follows:

\begin{prop}
  Let $(A,\sigma)$ be a central simple algebra of degree~$8$ with an
  orthogonal involution of trivial discriminant
  satisfying~\eqref{eq:C}, and let $W$ be a decomposition group of
  $(A,\sigma)$. The image $\overline{e_3(\sigma)}^W$ of
  $e_3(\sigma)\in H^3(F)/F^\times\cdot V$ in $H^3(F)/F^\times\cdot W$
  generates 
  $\ch_W$, and $f_3(\sigma)=f_3(W)$.
\end{prop}

\begin{proof}
  As above, let $Q$ be the quaternion division algebra
  Brauer-equivalent to $A$, so we may identify $A$ with $M_4(Q)$. Let
  $\rho$ be an involution on $M_6(Q)$ such that
  \[
  (M_6(Q),\rho)\in (A,\sigma)\orthsum
  \bigl((Q^+,\ba)\otimes(Q^-,\ba)\bigr).
  \]
  By definition, we have $e_3(\sigma)=\overline{e_3(\rho)}^V$ and
  $f_3(\sigma)=f_3(\rho)$. Now, 
  consider a decomposition of $(A,\sigma)$ with decomposition group
  $W$ (which necessarily contains $V$):
  \[
  (A,\sigma)\in \bigl((C_1^+,\ba)\otimes(C_1^-,\ba)\bigr) \orthsum
  \bigl((C_2^+,\ba)\otimes(C_2^-,\ba)\bigr).
  \]
  We have
  \[
  (M_6(Q),\rho) \in \bigl((C_1^+,\ba)\otimes(C_1^-,\ba)\bigr) \orthsum
  \bigl((C_2^+,\ba)\otimes(C_2^-,\ba)\bigr) \orthsum
  \bigl((Q^+,\ba)\otimes(Q^-,\ba)\bigr),
  \]
  which is a decomposition of $(M_6(Q),\rho)$ with decomposition
  group~$W$. Therefore, Theorem~\ref{hom.thm} shows that
  $\overline{e_3(\rho)}^W$ generates $\ch_W$ and
  $f_3(\rho)=f_3(W)$. The proposition follows because
  $f_3(\sigma)=f_3(\rho)$ and
  $\overline{e_3(\sigma)}^W=\overline{e_3(\rho)}^W$ since $V\subset W$.
\end{proof}

\Addresses

\end{document}